\documentclass[a4paper,10pt]{scrartcl}
\usepackage{geometry}
\usepackage[T1]{fontenc}
\usepackage{ae}
\usepackage[ngerman,english]{babel}
\usepackage[utf8x]{inputenc}
\usepackage{amsmath}
\usepackage{amssymb}
\usepackage{amsfonts}
\usepackage{dsfont}
\usepackage{array}
\usepackage{graphicx}
\usepackage{graphics}
\usepackage{algorithm}
\usepackage{algpseudocode}
\usepackage[caption=false]{subfig}
\usepackage[binary-units]{siunitx}
\usepackage{float}
\usepackage[hang,small,bf]{caption}
\usepackage{enumitem}
\usepackage{multirow}
\usepackage{fancyhdr}
\usepackage{colortbl}
\usepackage{booktabs}
\usepackage[e]{esvect}
\usepackage{xcolor}
\usepackage{eurosym}
\usepackage{textcomp}
\definecolor{darkblue}{rgb}{0,0,0.6}
\definecolor{revred}{rgb}{0.78,0,0}
\definecolor{darkgreen}{rgb}{0,0.6,0}
\definecolor{maroon}{rgb}{0.6,0,0}
\usepackage[hyperfootnotes=false]{hyperref}
\usepackage[textwidth=2.450cm, textsize=tiny]{todonotes}
\hypersetup{colorlinks=true, breaklinks=true, linkcolor=maroon, urlcolor=darkblue,citecolor=maroon}

\DeclareMathOperator*{\mydef}{\mathrel{\mathop:}=}

\DeclareMathOperator{\Exp}{\mathbb{E}}
\DeclareMathOperator{\Cov}{\mathbb{C}\text{ov}}

\DeclareMathOperator{\R}{\mathbb{R}}

\newcommand{\norm}[1]{\Vert #1 \Vert}
\newcommand{\sqnorm}[1]{\Vert #1 \Vert_2^2}

\newcommand{\relleetwo}{rel \, \ell_2}

\newcommand{\cD}{\mathcal{D}}
\newcommand{\cE}{\mathcal{E}}
\newcommand{\cF}{\mathcal{F}}

\newcommand{\cJ}{\mathcal{J}}

\newcommand{\cN}{\mathcal{N}}

\newcommand{\cS}{\mathcal{S}}
\newcommand{\cU}{\mathcal{U}}

\def\kwave/{\textbf{k}-\texttt{Wave}}

\newcommand{\termabb}[2]{\textit{#1} (\textit{#2})}

\title{High Resolution 3D Ultrasonic Breast Imaging by Time-Domain Full Waveform Inversion}
\author{Felix Lucka$^{1,2}$, Mailyn P\'{e}rez-Liva$^3$, Bradley E. Treeby$^4$, Ben T. Cox$^4$}

\date{%
    {\normalsize $^1$Centrum Wiskunde \& Informatica, Science Park 123, Amsterdam, The Netherlands\\[2pt]
    $^2$Centre for Medical Image Computing, University College of London, WC1E 6BT London, UK\\[2pt]
    $^3$Universit\'{e} de Paris, PARCC, INSERM, 56 Rue Leblanc, 75015, Paris, France\\[2pt]
    $^4$Department of Medical Physics and Biomedical Engineering, University College of London, WC1E 6BT London, UK\\[2pt]}
}

\begin{document}
\maketitle
\begin{abstract}
Ultrasound tomography (UST) scanners allow quantitative images of the human breast's acoustic properties to be derived with potential applications in screening, diagnosis and therapy planning. Time domain full waveform inversion (TD-FWI) is a promising UST image formation technique that fits the parameter fields of a wave physics model by gradient-based optimization. For high resolution 3D UST, it holds three key challenges: Firstly, its central building block, the computation of the gradient for a single US measurement, has a restrictively large memory footprint. Secondly, this building block needs to be computed for each of the $10^3-10^4$ measurements, resulting in a massive parallel computation usually performed on large computational clusters for days. Lastly, the structure of the underlying optimization problem may result in slow progression of the solver and convergence to a local minimum. In this work, we design and evaluate a comprehensive computational strategy to overcome these challenges: Firstly, we exploit a gradient computation based on time reversal that dramatically reduces the memory footprint at the expense of one additional wave simulation per source. Secondly, we break the dependence on the number of measurements by using source encoding (SE) to compute stochastic gradient estimates. Also we describe a more accurate, TD-specific SE technique with a finer variance control and use a state-of-the-art stochastic LBFGS method. Lastly, we design an efficient TD multi-grid scheme together with preconditioning to speed up the convergence while avoiding local minima. All components are evaluated in extensive numerical proof-of-concept studies simulating a bowl-shaped 3D UST breast scanner prototype. Finally, we demonstrate that their combination allows us to obtain an accurate $442\times442\times222$ voxel image with a resolution of \SI{0.5}{\milli\meter} using \textsc{Matlab} on a single GPU within 24 hours.
\end{abstract}

\newpage
\section{Introduction} \label{sec:Intro}

Screening using X-ray-based mammography has saved many lives by catching early-stage breast cancer. Nevertheless, there remain concerns related to tissue superposition, overdiagnosis, the effect of ionising radiation, the lower sensitivity in dense breasts, and the pain caused by breast compression \cite{MaWi2013}.
Alternative breast imaging approaches have therefore been explored, including MRI, which is expensive and time-intensive \cite{HeMo2019}, optical approaches such as photoacoustic tomography \cite{MaDa2019} and diffuse optical tomography \cite{ZhPo2020} which are limited by the significant scatting of light by tissue, and ultrasound (US). In conventional US imaging, high frequency acoustic waves are transmitted into the breast from a hand-held linear array of typically 128 detection elements. The back-scattered waves are detected by the same array and used to form qualitative 2D images in real time. It is often used as an adjunct modality to supplement mammography \cite{SoHa2019,GuFe2018}, but the operator-dependence means that hand-held US requires an expert user. It is also time-consuming when imaging a whole breast, making it unsuitable as a screening modality. In an attempt to overcome these drawbacks, automated breast US (ABUS) has been introduced, which uses a machine to scan a linear array perpendicular to its length over the breast to form a qualitative 3D reflection image \cite{Vo2019} but it does not give a quantitative image of tissue properties. In contrast, a principal aim in \textit{ultrasound tomography} (UST) of the breast \cite{CaMeScOu1981,DuLi18} is to obtain accurate high resolution \textit{quantitative} images of the tissue's acoustic properties. To achieve this, a pulse of US is transmitted into a pendant breast and measurements are made of the resulting time-varying acoustic field at multiple detector positions around the breast. This is repeated for multiple source positions (or, more generally, source distributions). The set of acoustic time series thus measured will incorporate both the unscattered (directly transmitted) and scattered parts of the acoustic field. To facilitate practical measurement times, large arrays are typically used, consisting of many hundreds of elements. In recent years, a number of such systems have been developed for UST breast imaging, including ones based on ring arrays \cite{DuLiScLiRo2013,DelphinusWebsite,OpPrSzSzSz2018}, rotating planar or linear arrays \cite{MaTeWiLe2018,QTUltrasoundWebsite,HuShCh2016}, and bowl-shaped detector arrays \cite{GeHoZaKaRu2017,PammothWebsite}. Once the data has been measured, the challenge is to use it to form accurate quantitative images of the tissue's acoustic properties. 

\subsection{Image Reconstruction for Ultrasound Tomography}
Typically three types of images have been produced from UST data: quantitative images of the acoustic attenuation and the speed of sound, and qualitative images indicating the `strength' of the scattering at each point in the image \cite{Co2007}. To obtain these images requires a nonlinear inverse scattering problem to be tackled. The various approaches to solving this, which can be characterised by the approximations and assumptions they make, represent different trade-offs between computational efficiency and imaging accuracy. UST systems and measurement protocols are therefore often designed and optimized for a particular image formation method, unlike other imaging modalities where it is typically the other way around. There is an extensive literature on UST reconstruction, and it shares much common ground with other areas of wave-based imaging, such as in seismology and non-destructive testing. This is not the place for a comprehensive account of that history, so we give just a few examples of typical reconstruction approaches, focusing on quantitative reconstructions of sound speed as that is the topic of this paper.

\textit{Travel-time} or \textit{time-of-flight} tomography, widely used in geophysics \cite{Ce2001}, assumes the wave propagation can be accurately described in terms of rays: either straight rays \cite{Sc2002,MaWiVaSlVa2007,JiPeRuFoDaZaJa2012,BiDaRuBe2014,OpPrSzSzSz2018}, thereby neglecting scattering, diffraction and refraction, or bent rays \cite{JoGrSaDuSj1975,An1990,LiDuLiHu2009, LiJaDiVaStMu2010,DaGeHaRu2012,JaLuCo2020}, which can account for refraction to some extent. The measured travel times are linearly related to the integral of the slowness (reciprocal of the sound speed) along the rays, and the system matrix, once constructed, can be inverted using standard linear algebra approaches. Ray-based approaches can lead to robust and computationally efficient algorithms to recover smoothed sound speed images, but usually fine details are lost. Because the significance of diffraction depends on the relative sizes of the scatterer and wavelength, the ray approximation becomes more accurate the higher the frequency used. However, at higher frequencies the sources and detectors become more directional, reducing the possible ray paths across the imaged region, and the acoustic attenuation of biological tissues is greater.

Another approach, sometimes called \textit{diffraction tomography}, considers the material heterogeneities to be weakly scattering perturbations over a background and uses Born or Rytov-type linearisations of the nonlinear scattering problem to arrive at linear methods of image formation. The free-space Green's function is used to estimate the effect of the perturbations on the field. In general for breast tissue, the weak scattering approximation does not hold, so adaptations of these linearisations are required to develop techniques suitable for quantitative breast imaging \cite{Ma2002,HuSi2011}. 

To solve the full nonlinear inverse scattering tomography problem requires an iterative approach. One way is to iteratively update both the forward operator and the unknown parameter in the linearised approximations above. In \textit{bent-ray reconstructions}, for example, both the rays and the unknown sound speed are updated at each iteration \cite{JaLuCo2020}. In the \textit{Distorted Born Iterative Method} \cite{HaEb1998,CaYuAl2019}, the Green's function and the unknown parameter are updated. This works well for low contrast media, but suffers from the failure of the Born series to converge for large contrast heterogeneities. 

Before we introduce the particular approach we follow in this work in more detail, we remark that as in many fields, there has been a recent flurry of interest into the question of whether deep neural networks can be utilized for ultrasonic imaging \cite{SlCoEl20}, see, e.g., \cite{Co1989,AlVeBuBoSa2018,FaYi2019,FeFrAn2019,LyTaDo2020,FaWaGeHe2020,FeMaFrAn2020,FeZwFrAn2020,JuBiDuScMa2020,KoDeDo2020,Zh2020,SlCoEl20} for approaches to improve the speed and accuracy of 2D UST reconstructions. While it seems likely that neural networks will find wide application in UST, there are reasons to think that, as in other imaging modalities, they will complement rather than supersede existing approaches \cite{ArMaOkSc19}. In particular in high resolution 3D settings, neural network components that link the object to be imaged to the measured data can hardly be fully learned. Efficient and accurate implementations will always contain components derived from physical models accurately describing wave propagation in biological tissue according to the best of our understanding.

\subsection{Full Waveform Inversion}

An approach that has gained increasing attention in recent years has been to formulate the image formation as a optimization problem, essentially fitting the measured data to a forward model of the acoustic propagation. This way of tackling the inverse scattering tomography problem is known as \textit{waveform tomography}, \textit{model-based inversion} or \textit{full waveform inversion} (FWI). Because the framework of FWI allows all the measured data to be used in the inversion - the whole time series, or all the frequencies - it can in principle produce quantitatively correct, high resolution images of the acoustic properties with many fewer measurements than is required for travel-time tomography (which essentially reduces each time series to a single value). There is a sizable literature on FWI in the seismic community, see \cite{Tr20} for a recent review. 

FWI constitutes a large-scale, non-convex, PDE-constrained optimization problem. In addition to the intrinsic difficulties of such optimization problems, iterative optimization schemes need to run numerous numerical wave simulations to compare artificial UST measurements caused by the current guess of the acoustic parameters to the real US data obtained. Depending on the scan design and the spatial resolution desired, this can present a very significant computational challenge. 

Dedicated hardware design can circumvent this bottleneck to some extent. For instance, when using a 2D ring-array focused in the plane, 2D FWI can be used to reconstruct the 3D volume slice-by-slice \cite{LiDuLiHu2009,DuLiScLiRo2013,GoRo2013}, but the risk of out-of-plane artefacts and the poor image resolution perpendicular to the plane makes this solution less than ideal \cite{SaWeLiRoDu17}. The approach taken by Wiskin et al. \cite{WiBoJoBeAbHa2007,WiBoJo2011,WiBoJoBe2012,WiBoIuKlLe2017} has been to design the hardware with a planar geometry so that a computationally-efficient approximation to the wave simulation can be exploited. 

There will be a highest usable frequency dictated by inherent attenuation through the breast and the sensitivity of the hardware, and this will limit the achievable image resolution. If this frequency is too high for the computations to be practicable then reducing the highest frequency, at the loss of the achievable image resolution, may be a compromise it is necessary to make. Goncharsky et al. \cite{GoRo2013,GoRoSe2016,GoRoSe2018,GoRoSe19b,GoRoSe19,GoSe20} follow this idea and consider using low-frequency transducers. They further combine with a ``2.5''D slice-by-slice reconstruction approach. 

However, there is not always the freedom to choose the hardware. For instance, in the case we are interested in, the scanner has been designed to be optimal for \textit{photoacoustic tomography} (PAT) of the breast \cite{PammothWebsite} and uses a hemispherical array with unfocused, broadband US transducers. UST, which will be an adjunct modality, will be performed using the same array and we aim to obtain the same, isotropic sub-millimeter resolution across the whole breast volume. 
Another example of a fully 3D UST system, also a bowl array, is the KIT 3D USCT system \cite{GeHoZaKaRu2017,GeBeHoZaTaBlLePeRu18}, which was designed for obtaining both high-quality reflectivity images and quantitative tissue images. To realize FWI for such UST scenarios, large high-performance computing clusters are typically used \cite{GoRo2013,BoMaViBaFi18,GoRoSe19} and/or the computation takes days, which severely limits the range of clinical applications for which these methods are viable. A notable exception is given by the work presented in \cite{BaTr20}, which we were made aware of during the writing of this article.

\subsection{Paper Scope and Structure} \label{subsec:Structure}
The key focus and contribution of our work is to develop and demonstrate a comprehensive computational strategy that achieves accurate high-resolution, 3D FWI for breast UST with a hemispherical array using only moderate computational resources (at least one GPU) and a self-imposed time limit of 24 hours. This would allow the method to fit into clinical trajectories for diagnosis and therapy planing at competitive costs compared to other modalities such as MRI.

In Section \ref{sec:Background} we summarize the mathematical modeling of UST and FWI, and describe the state of the art. Section \ref{sec:Methods} illustrates several separate computational challenges of FWI for high-resolution 3D UST, and describes how novel and established solutions can be combined into a comprehensive inversion strategy. Extensive numerical studies in Section \ref{sec:NumericalStudies} first validate solutions to the separate sub-problems before the results of the whole scheme are demonstrated in Section \ref{sec:FullDemo}. Finally, we discuss the results of our work and point to future directions of research in Section \ref{sec:DisOut} before closing with conclusions in Section \ref{sec:Con}. Table \ref{tbl:Abb} lists all commonly occurring abbreviations for reference.

\begin{table}[h!]
\caption{\label{tbl:Abb} List of commonly occurring abbreviations.} 
\begin{tabular}{@{}lll}        
\toprule
Abbreviation & Meaning & Reference  \\
\midrule
FD       & \textit{Frequency Domain}            & Sec. \ref{subsec:FWI} \\
FWI      & \textit{Full Waveform Inversion}     & Sec. \ref{subsec:FWI} \\
US       & \textit{Ultrasound}                  & Sec. \ref{sec:Intro}\\
UST      & \textit{Ultrasound Tomography}       & Sec. \ref{subsec:UST}\\
SGD      & \textit{Stochastic Gradient Descent} & Sec. \ref{subsec:SGDandSE} \\
SLBFGS   & \textit{Stochastic L-BFGS Method}    & Sec. \ref{subsec:SGDandSE} \\
TD       & \textit{Time Domain}                 & Sec. \ref{subsec:FWI} \\
TR       & \textit{Time Reversal}               & Sec. \ref{subsec:TRGradient} \\
\bottomrule
\end{tabular}
\end{table}

\section{Full Waveform Inversion for Ultrasound Tomography} \label{sec:Background}

This section presents the model used to describe ultrasound propagation in soft tissue (Sec.\ \ref{subsec:US}) and the mathematical formulation of the inverse problem of UST (Sec.\ \ref{subsec:UST}). The FWI approach to tackling this problem, the progress that has been made, and the challenges that remain, are then described (Secs.\ \ref{subsec:FWI} and \ref{sec:FWI-Challenges}).

\subsection{Modelling Ultrasound Propagation in Soft Biological Tissues} \label{subsec:US}

It is common to model soft biological tissues as compressible fluids and describe acoustic propagation therein by linearising the equations of fluid dynamics derived from conservation laws and combining them with an empirical equation of state. Following this approach, a sufficiently low amplitude ultrasonic wave can be modelled by the following system of first order equations \cite{TrJaReCo12}:
\begin{align}
\frac{\partial}{\partial t} v(x,t)& = - \frac{1}{\rho_0(x)} \nabla p(x,t) \hspace{18.2em}  \text{\small (momentum conservation)}\label{eq:MonCon} \\
\frac{\partial}{\partial t} \rho(x,t)& = - \rho_0(x) \nabla \cdot v(x,t) - v(x,t) \cdot \nabla \rho_0(x) + m(x,t)  \hspace{5.8em}  \text{\small (mass conservation)} \label{eq:MassCon}\\ 
p(x,t)& = c_0^2(x) (\rho(x,t) + d(x,t)\cdot \nabla\rho_0(x) - L \rho(x,t)),  \hspace{6.7em}   \text{\small (pressure-density relation)} \label{eq:PreDenRel}
\end{align}
where $v$ is the acoustic particle velocity (the time derivative of the acoustic displacement $d$), $p$ and $\rho$ are the acoustic pressure and acoustic density fluctuations, $\rho_0$ and $c_0$ are the ambient density and sound speed, and $m$ is a mass source term. $L$ models acoustic absorption and dispersion using fractional Laplacians:
\begin{align}
L &= \tau(x) \frac{\partial}{\partial t} \left( - \nabla^2  \right)^{\tfrac{y}{2}-1} + \eta(x) \left( - \nabla^2  \right)^{\tfrac{y+1}{2}-1},  \\
\tau(x) &= - 2 \alpha_0(x) c_0(x)^{y-1}, \hspace{6em} \eta(x) = 2 \alpha_0(x) c_0(x)^y \tan(\pi y/ 2),
\end{align}
where $\tau$ and $\eta$ denote absorption and dispersion proportionality coefficients, $\alpha_0$ is the power law prefactor, and $y$ is the power law exponent \cite{TrCo11b}. Equations \eqref{eq:MonCon}-\eqref{eq:PreDenRel} can be combined into the following lossy second-order wave equation for a heterogeneous medium:
\begin{equation}
A  \, p(x,t) \mydef \left( \frac{1}{c_0^2(x)}\frac{\partial^2}{\partial t^2}  - \rho_0(x)\nabla\cdot \left(\frac{1}{\rho_0(x)} \nabla \right) + L \nabla^2 \right) p(x,t) =  \frac{\partial}{\partial t} m(x,t) =: s(x,t), \label{eq:WaveSecondOrder}
\end{equation}
where we introduced a lossy d'Alembert operator $A$. As initial conditions we have $p = 0$ and $\partial_t p  = 0$. Suitable data pre-processing typically allows us to model the propagation as if it occurs in an unbounded domain, i.e., no explicit boundary conditions are required (reflections from experimental equipment can be time-gated out and waves propagating deeper into the body are absorbed). See \cite{TrCo10,TrJaReCo12,ArBeCoLuTr16} for more detailed discussions on the mathematical modeling and  \cite{BeMoKoLa17,GuCaTaNaWa20,TaSoNaJaVeDo20} for recent extensions of FWI to UST involving hard tissues such as bone.

\subsection{The Inverse Problem of Ultrasound Tomography} \label{subsec:UST}

In most UST systems, movable arrays of US transducers are placed around the sample $\Gamma \subset \R^3$ and acoustically coupled with it. For a single recording, a subset of transducers emit an acoustic pressure signal of length $T_s$ while the others are receiving for a time $T > T_s$, chosen long enough such that the acoustic pressure remaining inside the scanner is not measurable anymore. A whole scan consists of $n_s$ such recordings, between which the transducer arrays may be moved. Here, we model a UST scan in the following way: we have $n_s$ temporal sources $s_i(x, t)$, with $\mathrm{supp}(s_i) \subset \Gamma^c \times [0,T_s]$ (in practice, each source is only defined over a small region corresponding to the front face of the transducer). The measurement process is modelled by applying a linear operator $M_i$ to the pressure fields $p_i(x,t), t\in [0,T]$ resulting from \eqref{eq:WaveSecondOrder}, i.e., we have
\begin{equation}
A(u) \, p_i(x,t) = s_i(x, t),  \hspace{3em} f_i = M_i p_i, \hspace{3em} i=1, \ldots,n_s. \label{eq:USTmodel}
\end{equation}
Here, $u$ denotes the unknown acoustic material properties inside the sample $\Gamma$, in the extreme case $u = (c_0, \rho_0, \alpha_0, y)$, where $c_0$, $\rho_0$ and $\alpha_0$ depend on $x$. We assume that all other properties are sufficiently well-known to ignore their modeling error. The general inverse problem of UST is to recover $u$ (or features of it) given a noisy data set $\{f^\delta_i\}_{i=1}^{n_s}$. A general discussion of the uniqueness and stability of this problem can be found in \cite{Ta05}. As mentioned in the Introduction, many different approaches have been used to tackle this problem.

\subsection{Full Waveform Inversion} \label{subsec:FWI}

In FWI, we assume that for a given $u$, we can solve \eqref{eq:USTmodel} to simulate data, i.e., $f_i(u) := M_i A(u)^{-1} s_i$. Then, we try to optimize $u$ such that the discrepancies between simulated and measured data become small:
\begin{equation}
\min_{u \in \cU} \, \cJ(u) , \qquad \text{with} \qquad \cJ(u) := \sum_i^{n_s} \cD_i(u)  :=  \sum_i^{n_s} \cD \left( f_i(u) , f^{\delta}_i \right) =  \sum_i^{n_s} \cD \left( M_i A^{-1}(u) s_i  , f^{\delta}_i \right),  \label{eq:FWI}
\end{equation}
where $\cD(f, g)$ is a loss function (see \cite{Ya18,EnRnYa20} for a discussion of suitable loss functions), and $\cU$ is a set of constraints on $u$, e.g., bound constraints. First-order optimization schemes solve \eqref{eq:FWI} using only the gradient $\nabla \cJ(u)$, which is given by the sum over terms of the form $\nabla_u \cD \left( M A^{-1}(u) s  , f^\delta \right)$. Such expressions can be computed efficiently using the adjoint state method \cite{Pl06}, which we will summarize here in a short, instructive way. One starts by differentiating both sides of \eqref{eq:USTmodel} wrt $u$:
\begin{align}
&\frac{\partial A}{\partial u} p + A \frac{\partial p}{\partial u} = 0 \quad \Rightarrow \quad \frac{\partial p}{\partial u} = - A^{-1}  \frac{\partial A}{\partial u} p \quad \Rightarrow \quad \frac{\partial f}{\partial u} = - M A^{-1}  \frac{\partial A}{\partial u} p \nonumber \\
&\Rightarrow \quad \frac{\partial \cD}{\partial u} = \left(\frac{\partial f}{\partial u} \right)^T \frac{\partial \cD}{\partial f} = - \left( \frac{\partial A}{\partial u} p \right)^T A^{-T} M^T  \frac{\partial \cD}{\partial f}
\end{align}
Here, $A^{-T}$ is the operator solving the adjoint wave equation \cite{Pl06}. For instance, for $u = c_0$ and $L = 0$, we get with $\frac{\partial A}{\partial c_0} = - \frac{2}{c_0^3} \frac{\partial^2}{\partial t^2}$:
\begin{align}
&\nabla_{c_0} \cD \left( f(c_0) , f^{\delta} \right) =  \int_0^T \frac{2}{c_0(x)^3} \left(\frac{\partial^2 p(x,t)}{\partial t^2} \right) q^*(x, T - t) \, dt, \label{eq:Gradc0}
\end{align}
where $q^*$ is the adjoint pressure field obtained by solving $A(u) q^* = s^*$, and $s^*(x, t)$ is the time-reversed data discrepancy gradient $\frac{\partial \cD}{\partial f}$ \cite{Pl06}. We will use $\cD(f,g) = \tfrac{1}{2} \sqnorm{f-g}$, for which we have $\frac{\partial \cD}{\partial f} = f-g$. Corresponding formulas for $\rho_0$, $\alpha_0$ and $y$ are listed in Appendix \ref{sec:DenAbsGrad}. 

The focus of this work is the practical feasibility of computing a sufficiently accurate and spatially resolved approximation to the solution of the core problem \eqref{eq:FWI} for 3D breast UST within a reasonable amount of time and with reasonable computing resources. For this reason, we will only include bound constraints on $u$ via $\cU$ in the main studies. To embed more sophisticated \textit{a-priori} knowledge on $u$ and account for noise models and model discrepancies in $f^\delta$, one needs to extend \eqref{eq:FWI} by adding regularization functionals (e.g., to penalize unwanted spatial features of $u$) or additional constraints (see \cite{MaWaLiDuAn17,EsGuLeArHe18} and references therein). The additional noise stability study in Appendix \ref{sec:Noise} will demonstrate the use of a simple regularization strategy.

We introduced FWI in the \termabb{time domain}{TD} here, which is typically also the domain the measurements are obtained in. One can also formulate FWI in the \termabb{frequency domain}{FD}. Then, computing $\cJ(u)$ and its gradient requires solving a direct and an adjoint Helmholtz equation for each frequency. If the emitted waves are narrow band and/or the transducers are only sensitive within a narrow band (modelled by $M_i$), this can be very advantageous but typically leads to high memory requirements if high resolution in 3D is desired (see \cite{SaWeLiRoDu17} and references therein). 
 
Depending on the scanning geometry, transducer characteristics, and measurement protocol and whether TD-FWI or FD-FWI is pursued, different numerical schemes are advantageous for solving the wave or Helmholtz equation. Direct methods such as finite-difference, pseudospectral, finite/spectral element methods or discontinous Galerkin methods aim to solve the PDEs directly while integral equation methods such as boundary element methods first transform them. Asymptotic or approximate methods such as geometrical optics or Gaussian beams only capture a part of the wave-matter interaction while typically being computationally much more efficient. See \cite{Fi11,Ig16} for overviews on computational wave propagation methods for seismic FWI imaging.

\subsection{Challenges of High Resolution 3D Time Domain Full Waveform Inversion} \label{sec:FWI-Challenges}

Many UST systems are designed to capture 3D information and emit unfocused, broadband waves to obtain sub-millimeter isotropic spatial resolutions over the whole breast volume \cite{GeHoZaKaRu2017,PammothWebsite}. For instance, within the PAMMOTH project \cite{PammothWebsite}, we built a scanner with a bowl-shaped transducer array and aim to reach a resolution of $\leq$\SI{0.5}{\milli\meter}. Even when using the adjoint state method, solving TD-FWI \eqref{eq:FWI} for such scenarios is computationally challenging: For the PAMMOTH setup, even an efficient \textit{k}-space pseudospectral method \cite{TrCo10} takes at least 10 minutes to solve a single wave simulation on a recent GPU. As $\cJ(u)$ involves the sum over $n_s$ sources \eqref{eq:FWI}, computing $\nabla \cJ(u)$ requires $n_s$ parallel gradient computations \eqref{eq:Gradc0} each involving two wave simulations. On a single GPU, this would take $14$ days for $n_s = 1024$. While using a cluster with many computational nodes could reduce this \cite{BoMaViBaFi18,GoRoSe19}, a closer look at \eqref{eq:Gradc0} reveals a problem of the TD adjoint state method: Ideally, the field $p(x,t)$ needs to be kept in GPU memory or at least in the main memory of the computational node while the adjoint field $q^*(x,t)$ is computed. For the PAMMOTH setup, this would require up to \SI{150}{\giga\byte} working memory. Current GPU cards are limited to \SI{48}{\giga\byte} and typical nodes in computational clusters are not equipped with that much main memory, either. In addition to these difficulties, first-order optimization methods for FWI \eqref{eq:FWI} suffer from slow convergence and may get stuck in suboptimal local minima of the non-convex function $\cJ(u)$. 

In the following section, we describe techniques to circumvent each of these problems and combine them into a comprehensive computational strategy to achieve high resolution 3D TD-FWI for ultrasonic breast imaging.

\section{Improving Memory Footprint, Efficiency and Convergence} \label{sec:Methods}

\subsection{Memory-Efficient Gradient Computation Using Time Reversal} \label{subsec:TRGradient}

First, we describe an approach to circumvent storing $p(x,t)$, $x \in \Gamma$, $t \in [0,T]$. Note that in UST, all sources and sensors are outside $\Gamma$ and $p(x,0) = \partial_t p(x,0) = 0$. We first consider the case of no absorption, i.e., $L = 0$. By Huygens' principle, the field inside $\Gamma$ is uniquely determined by its trace on the boundary $\partial \Gamma \times [0,\infty]$. The theory of \termabb{time reversal}{TR} describes under which conditions the time reversal symmetry of the wave equation \eqref{eq:WaveSecondOrder} allows $p(x,t)$ to be reconstructed by simulating the following wave equation for the time-reversed pressure $p^{\triangleleft}(x,t)$: 
\begin{equation}
A(u) p^{\triangleleft}(x,t) = 0 \hspace{1em} \text{in} \hspace{.5em}\Gamma \times [0,T^{\triangleleft}], \hspace{1em}
p^{\triangleleft}(x,t)   = p(x,T^{\triangleleft}-t) \hspace{1em} \text{on} \hspace{.5em}\partial\Gamma \times [0,T], 
\hspace{1em} p^{\triangleleft}(x,0) = \partial_t p^{\triangleleft}(x,0) = 0 \hspace{1em} \text{in} \hspace{.5em} \Gamma \label{eq:TRBwdWave}
\end{equation}
Namely, $p^{\triangleleft}(x,t)$ is a good approximation of $p(x,T^{\triangleleft} - t)$ if $\Gamma$ is convex, $u$ is \textit{non-trapping} and $T^{\triangleleft}$ is chosen large enough such that $p(x, T^{\triangleleft})$ and $\partial_t p(x, T^{\triangleleft})$ are small and monotonically decaying inside $\Gamma$ \cite{Fi92,WuFi92,CaFi92,KuKu11}. In 3D UST of the breast, these conditions are fulfilled: $\Gamma$ can be chosen as a convex set including the breast (e.g., in the PAMMOTH system, the breast is placed in a plastic cup with convex shape). The acoustic parameter variations are such that no waves get trapped in the sample and after all acoustic energy has entered the breast, it decays exponentially fast, see Appendix \ref{sec:SourceEncodingTimeInvariantSys}. Choosing $T^{\triangleleft} = T$ is a safe choice as we chose $T$ large enough to ensure the pressure inside the scanner is not measurable any more (cf. Section \ref{subsec:UST}). Note that \eqref{eq:TRBwdWave} is not the same as the adjoint wave equation \cite{ArBeCoLuTr16}.

TR was developed for focusing ultrasound waves through inhomogenous media \cite{Fi92,WuFi92,CaFi92} and is used for image reconstruction in \termabb{Photoacoustic Tomography}{PAT} \cite{KuKu11,Bea11}. Here, we use it as a numerical trick to approximately replay the forward field $p(x,t)$ backwards in time, in parallel to solving the adjoint wave equation. This way, we do not need to keep the pressure field $p(x,t)$ in the sample's volume $\Gamma$ in memory, only its trace on the boundary $\partial\Gamma$, which reduces the memory footprint considerably. Note that similar techniques have been developed in the seismic literature. See \cite{Mu18} and references therein for the idea to shift the storage onto the boundary for our type of boundary conditions.

The integral in \eqref{eq:Gradc0} can directly be accumulated during the parallel time stepping scheme solving the both TR and adjoint wave equations: The state variables in each computation are exactly what is needed to compute the contribution to the integral. For instance, the first time step of the TR wave simulation results in $p^{\triangleleft}(x,0)$, which is approximately $p(x,T)$, while the first time step of the adjoint wave simulation results in $q(x,0)$. With these variables, the contribution to the integral \eqref{eq:Gradc0} for time $t=T$ can be computed. The price we pay for the heavily reduced memory footprint is that we have to run three instead of two wave simulations to compute \eqref{eq:Gradc0}. However, the simulation domain of the TR wave simulation only needs to enclose $\Gamma$ and can therefore often be chosen smaller. Appendix \ref{sec:kSpaceImp} describes the implementation of the TR-based gradient computation using the \textit{k}-space pseudospectral method implemented in the \kwave/ Matlab Toolbox \cite{TrCo10} and Section \ref{subsec:ValTRGrad} validates it numerically. 

In the case of absorption, i.e., $L \neq 0$, \cite{TrZhCo10} describes how to modify TR to recover the pressure fields (essentially, the sign of the absorption needs to be switched). Note however, that compared to the application in PAT image reconstruction considered in \cite{TrZhCo10}, no additional regularization is needed in our case as the boundary time series are not contaminated with measurement noise.

\subsection{Tuneable Stochastic Gradient Optimization Using Delayed Source Encoding} \label{subsec:SGDandSE}

FWI \eqref{eq:FWI} is a \textit{finite sum minimization} problem. As such, one can approximate $\nabla \cJ(u)$ by
\begin{equation}
    g_{\cS}(u) :=  \frac{n_s}{|\cS|} \sum_{j \in \mathcal{S}} \nabla \mathcal{D}_j(u), \qquad \cS \subset \{1,\dots, n_s \}. \label{eq:SubSamGrad}
\end{equation} 
In \textit{incremental gradient} or \textit{ordered sub-set} methods \cite{Be15}, one chooses a different subset $\cS_k$ in each iteration $k$ in a predetermined way. If the subsets are chosen at random, $g_{\cS}(u)$ becomes a stochastic estimator of $\nabla \cJ(u)$ and \termabb{stochastic gradient descent}{SGD} schemes need to be used. Fortunately, supervised training of deep neural networks are also instances of finite sum minimization problems (\textit{empirical risk minimization}) and the recent success of deep learning techniques has stimulated research into efficient SGD schemes \cite{BoCuNo18}. In Section \ref{subsec:ValStochOpt}, we will demonstrate the benefits of using a recently developed \termabb{stochastic L-BFGS}{SLBFGS} method (see \cite{FaGlGi17} and Appendix \ref{sec:SLBFGS}) for UST over earlier schemes.  

Besides being unbiased, i.e., $\Exp \left[g(u)\right] = \nabla \cJ(u)$, a desirable property of a gradient estimator $g(u)$ is that it is both computationally efficient (fast to compute) and stochastically efficient (small variance/error). Computing $g_{\cS}(u)$ with a subset of size one takes only two wave simulations (or three if TR is used). It turns out that due to the structure of \eqref{eq:FWI}, there are gradient estimators with the same computational but higher stochastic efficiency. In this work, we assume $M_i := M$. Let $w_i$, $i=1,\ldots,n_s$ be random weights with $\Exp \left[ w \right] = 0$, $\Cov[w] = I$ and define the synthetic source activations $\hat{s} := \sum_i^{n_{s}} w_i s_i$ and corresponding synthetic data $\hat{f}^\delta := \sum_i^{n_{s}} w_i f^\delta_i$. Then, $g_{\cE}(u) :=  \nabla_{u} \cD \left( M A(u)^{-1} \hat{s} , \hat{f}^{\delta} \right)$ is an unbiased estimator of $\nabla \cJ(u)$. Practically, we run two ``super-wave'' simulations where all sources are fired simultaneously. This technique exploits the linearity of the wave equation and is called \textit{source encoding}, see \cite{HaChHe12} for the general theory, \cite{WaMaAnLiDuAn15} for the form we use it here (but we apply it in 3D), and Appendix \ref{sec:SourceEncodingTimeInvariantSys} for a more detailed description. Choosing $w_i = \pm 1$ with equal probability (\textit{Rademacher distribution}) minimizes the variance of $g_{\cE}$ and will be used from now on (note that one could also define $w$ such that $g_{\cS}$ is recovered). The case $M_i \neq M$ will be examined in forthcoming work that applies the techniques discussed here to concrete experimental data scenarios.

While SGD methods need to reduce the variance of the update progressively in order to converge at all, tailored variance reduction schemes can improve the convergence rate to that of deterministic methods (cf. Section 5 in \cite{BoCuNo18}). If the variance of the gradient estimator can not be influenced in a direct way, variance reduction can only be achieved via reducing the step size or averaging iterates or gradients. Having fine control over the variance of $g(u)$ is more advantageous. For $g_{\cS}(u)$, one can only increase the size of $\cS$; for $g_{\cE}$, one would average multiple realizations. In both cases, we would have a sudden increase of the computational effort which means a rather coarse control over the variance \cite{HeBoAfThKrTrFi20}. If multiple GPUs are available, one can perform these computations in parallel, which we will investigate in Section \ref{subsec:MultipleGPUs}. However, if different GPU architectures are used in this way, the slowest one limits the overall efficiency, cf. Section \ref{subsec:GPUsCmp}. To overcome these problems, we propose a novel source encoding scheme $g_{\cE}^\tau$ with fine variance control that also exploits the time-invariance of the wave equation and the fact that all waves eventually leave the region of interest (which is essential for the sequential scanning of UST, cf. Section \ref{subsec:UST}): Let $d_i \sim U[0,1]$ be random, $\tau \geq 0$ a fixed maximal delay and define 
\begin{equation}
    \hat{s}(x,t) := \sum_i^{n_{s}} w_i s_i(x, t - d_i \tau), \qquad \hat{f}^\delta(x,t) := \sum_i^{n_{s}} w_i f^\delta_i(x, t - d_i \tau) \label{eq:DelaySE}
\end{equation} 
For $\tau = 0$, we recover the conventional source encoding where all sources are fired simultaneously while for $\tau > 0$, we activate the sources with random delays. While this necessitates running the super-wave simulation for a longer time $T + \tau$, it dampens the cross-talk between the pressure fields coming from different sources, which is the origin of the error of conventional source encoding. For $n_{s} = 2$ and $\tau = 2T$, the average delay between the two sources is $T$ and we essentially fire sources with the delay we deemed necessary for regarding the measurements as completely separate (cf. Section \ref{subsec:UST}). More generally, $g_{\cE}^\tau$ converges to $\nabla \cJ(u)$ for growing $\tau$ for any realization of Rademacher weights $w$ and $d$. More details can be found in Appendix \ref{sec:SourceEncodingTimeInvariantSys}. Note that our approach extends the 2D FWI random time-delay approach described in \cite{ZhHu14} by combining it with the rigorous  stochastic analysis from \cite{HaChHe12}. Recently, another interesting TD approach to reduce the cross-talk between different sources by running longer ``super-wave'' simulations has been presented in \cite{TrBa19,BaTr20}: The encoding assigns a monochromatic time course to each source $i$, chosen from a set of regularly spaced frequencies. The forward and adjoint wave simulations are run long enough to produce the steady-state pressure fields from which the contributions of each source can be decoded by Fourier decomposition.

In Section \ref{subsec:ValStochOpt}, we will validate the proposed gradient estimator numerically. On a single GPU, this technique allows one to increase the gradient estimator's precision in arbitrarily small amounts while on heterogeneous multi-GPU environments, it allows faster GPUs to compute more precise gradient estimates using tailored $\tau > 0$ while the slowest GPUs use $\tau = 0$.

\subsection{Improving Convergence Using Coarse-To-Fine Schemes and Preconditioning} \label{subsec:MultGrid}

FWI \eqref{eq:FWI} is a non-convex optimization problem. Iterative optimization techniques based on local descent directions like gradient descent schemes may experience slow convergence or even get stuck in local minima that only explain parts of the data. In seismic imaging, the latter is known as \emph{cycle skipping} \cite{Ya18,EnRnYa20,Tr20}. A common way to avoid it in FD-FWI schemes is to restrict the frequency range to the lower frequencies first and and increase it progressively \cite{WiBoJoBeAbHa2007}. In TD schemes, this corresponds to solving the FWI on a coarse spatio-temporal grid first, which is then progressively refined. In the most basic case, the interpolation of the solution computed on the coarse grid is used to initialize the optimization scheme on the finer grid but more sophisticated multi-grid schemes can be derived \cite{JaHo17}. In the \textit{k}-space pseudospectral method used here, all grids are regular and implementing correct up and down-sampling is easy, see Appedix \ref{sec:MultiGridDetails}. In 3D, a coarsening of the spatio-temporal grid by factor of $\eta = 2$ (i.e., $dx^c = \eta \, dx^f$, $dt^c = \eta \, dt^f$,...) leads to a reduction of the computational operations by a factor of $\eta^4 = 16$. We will illustrate in Section \ref{subsec:ValMultGrid} that it is essential to exploit this advantage to obtain fast computational schemes. 

Preconditioning techniques try to improve the convergence of numerical optimization schemes by reformulating the original problem. For instance, instead of using $u = c_0^2(x)$ in \eqref{eq:FWI}, we could solve for $u = 1/c_0(x)$ (\textit{slowness}), $u = 1/c_0^2(x)$ (\textit{squared slowness}) or a weighted version thereof, $u = w(x)/c_0^2(x)$. The weights $w(x)$ could be a function of the average distance of $x$ to the transducer locations. Due to the non-convex nature of $\cJ(u)$ and the way SLBFGS constructs low-rank quadratic approximations to $\cJ(u)$ based on stochastic function and gradient evaluations, this choice is not trivial. In Section \ref{subsec:ValPrecon}, we will examine different possibilities.

\section{Numerical Accuracy and Efficiency Studies} \label{sec:NumericalStudies}

\subsection{Setup: 3D Digital Phantom and Transducer Arrangement}

\begin{figure}[tb!]
\centering
\hfill 
\subfloat[][\label{subfig:SetupA}]{\includegraphics[height=0.5\textwidth]{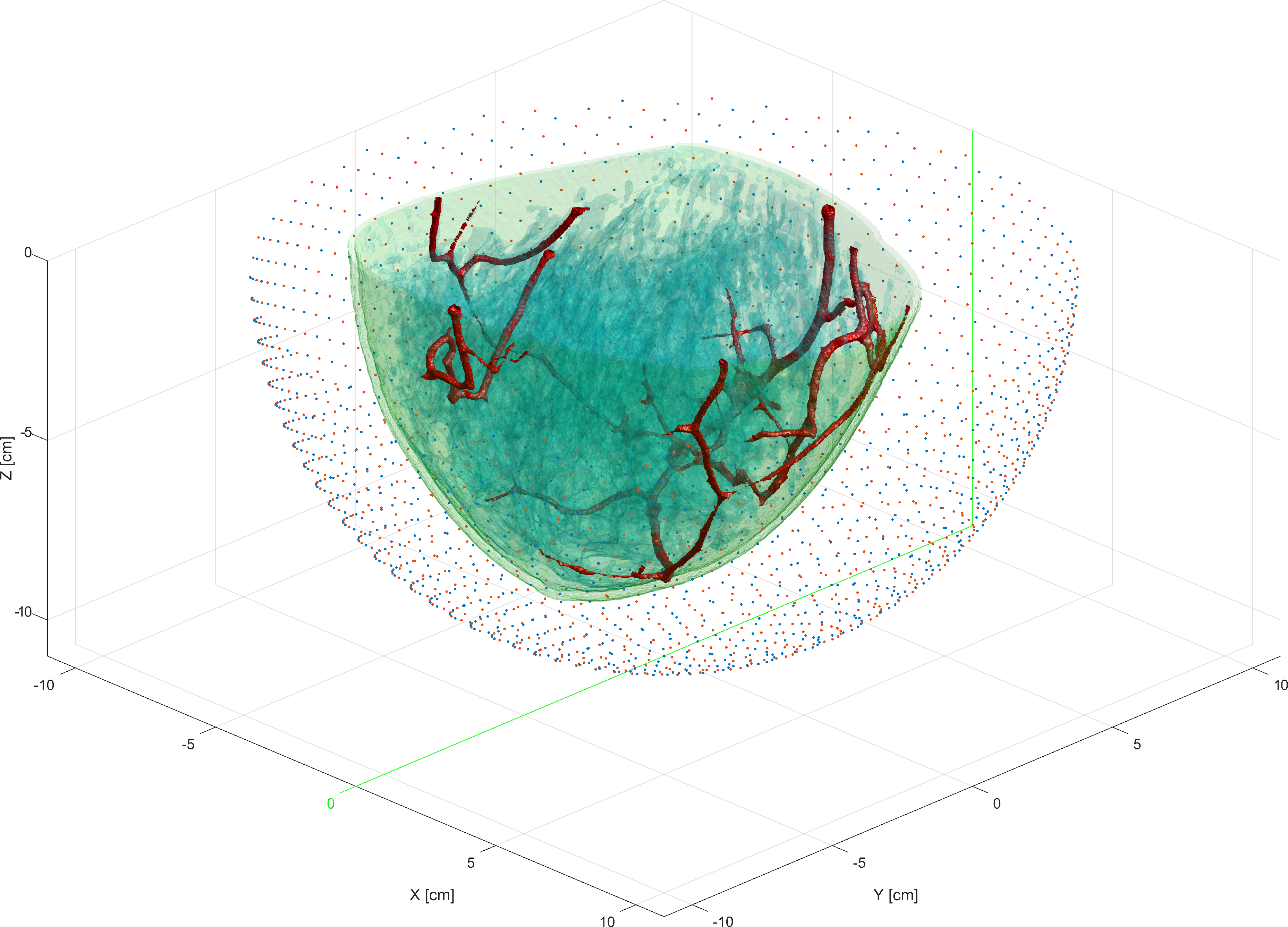}}
\hfill
\subfloat[][\label{subfig:SetupB} $x = 0$]{\includegraphics[height=0.5\textwidth]{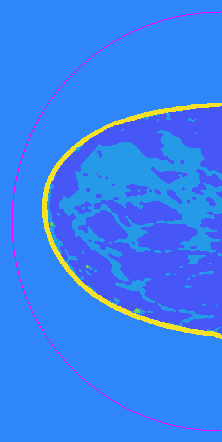} \hskip 2pt \includegraphics[height=0.5\textwidth]{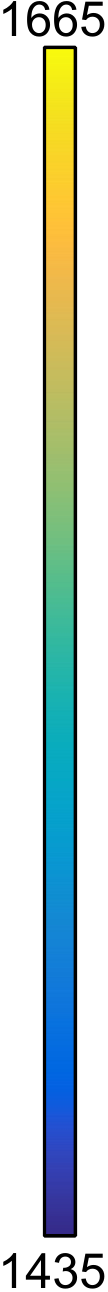}}
\hfill

\caption{\protect\subref{subfig:SetupA} Illustration of the setup used in the numerical studies: The surfaces of different breast tissues have been extracted from a high-res version of the voxel-grid-based breast phantom and rendered with different colors and transparencies. Shown here are blood vessels (red, in-transparent), fibro-glandular tissue (blue, transparent) and skin (yellow, transparent). Source and sensor positions are indicated by the blue and orange dots distributed on a half-sphere, respectively.  \protect\subref{subfig:SetupB} Illustration of the central x-slice of the 3D volume (corresponding to the green line in \protect\subref{subfig:SetupA}). Vertical and horizontal axis correspond to $y$ and $z$ axis, respectively. Different materials are colour coded depending on their speed of sound values ($c_0$), ranging from fat (\SI{1470}{\meter / \second}) over water/background (\SI{1500}{\meter / \second}), fibro-glandular tissue (\SI{1515}{\meter / \second}), blood vessels (\SI{1584}{\meter / \second}) to skin (\SI{1650}{\meter / \second}). The pink pixels indicate the half-sphere on which the transducers are located, the reconstruction region $\Gamma$ is composed by all pixels containing breast tissues.}
   \label{fig:Setup}
\end{figure}

The setup used in the numerical validation studies resembles the PAMMOTH scanner and is illustrated in Figure~\ref{fig:Setup}. The computational domain is $221\times 221 \times 111$\si{\milli\meter} wide and the ultrasound transducers are located on a half-spherical surface of radius \SI{104.5}{\milli\meter}. We want to compute a FWI reconstruction at spatial resolution $dx =$ \SI{0.5}{\milli\meter}, which leads to a computational domain of $442 \times 442 \times 222$ voxels. The speed of sound bounds enforced during the optimization are $c_0^{min} =$\SI{1350}{\meter / \second} and $c_0^{max} =$\SI{1800}{\meter / \second}.  Based on this and using a $cfl$ number of $0.3$, we set $\Delta t = \frac{cfl \cdot  \Delta x}{c_0^{max}} =  83.3$\si{\nano\second}. We set $T$ to twice the length of the computational domain divided by $c^{min}_0$, which leads to $T/\Delta t = 3912$ time steps. 

The anatomically realistic numerical breast phantom used is part of the OA-Breast Database\footnote{\href{https://anastasio.bioengineering.illinois.edu/downloadable-content/oa-breast-database/}{https://anastasio.bioengineering.illinois.edu/downloadable-content/oa-breast-database/}}, which was designed for (photo-)acoustic simulations and is described in \cite{LoZhMaApAn17}. It is based on the tissue segmentation of a clinical contrast enhanced MRI of a breast in prone, free-hanging position. The original segmentation is fitted and interpolated into our setup. The tissue-dependent speed of sound values we use are the same as in \cite{LoZhMaApAn17} and can be found in the caption of Figure~\ref{fig:Setup}. The mass density $\rho_0$ was assumed constant at \SI{1000}{\kilo\gram / \meter^3} for simplicity here. 

We distribute $2048$ transducer locations over the half-sphere by the golden section method. Half of them will be emitting (sources, $n_s = 1024$) while the other half will be receiving (sensors). Simulated data of this artificial source-sensor setup is an array of size $1024 \times 1024 \times 3912$ which corresponds to \SI{15.28}{\giga\byte} in single precision. 
We assume idealized spatio-temporal transducer characteristics here, i.e., we model them as point-like in space and their temporal impulse response is a delta function. The sources emit a single broadband pressure pulse (``click'') smoothed in time to remove temporal frequencies not supported by the computational grid (the grid supports up to \SI{1.5}{\mega\Hz} at $c_0 = $\SI{1500}{\meter / \second} \cite{TrCo10}). More details can be found in Appendix \ref{sec:MultiGridDetails}.

Note that the setup we use here was primarily designed to validate the accuracy and efficiency of the computational solutions we described in the previous section in the best possible way. For this reason, it reproduces the key properties of 3D FWI for breast imaging but omits modeling features of real-world systems that would complicate the interpretation of our results. It is furthermore not designed to examine the condition of the inverse problem \ref{eq:USTmodel}.

\subsection{Baseline Reconstructions} \label{subsec:BaselineRec}

\begin{figure}[htb!]
   \centering
\hfill
\subfloat[][\label{subfig:Ini}]{\fbox{\includegraphics[height=0.4\textwidth]{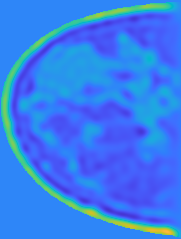}} \hskip 2pt \includegraphics[height=0.4\textwidth]{ParulaNew1435-1665.pdf}}
\hfill
\subfloat[][\label{subfig:IniDiff}]{\fbox{\includegraphics[height=0.4\textwidth]{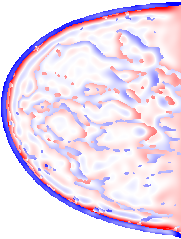}} \hskip 2pt \includegraphics[height=0.4\textwidth]{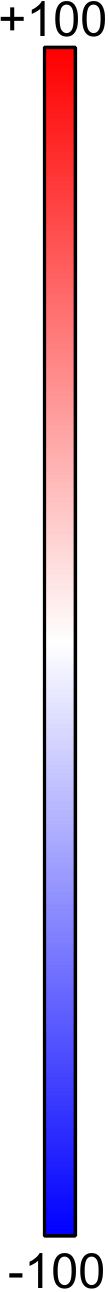}}
\hfill
\caption{\protect\subref{subfig:Ini} Zoom into central x-slice of the FWI solution computed using a grid size of \SI{2}{\milli\meter}, $c_0^{ini}$. Colormap in m/s. All reconstructions will be displayed in this way. \protect\subref{subfig:IniDiff} Error $c_0^{ini} - c_0^{\dagger}$ displayed with bidirectional colour map in m/s.}
   \label{fig:Ini}
\end{figure}

Before we can examine the different aspects discussed in Section \ref{sec:Methods} in depth, we need to illustrate some basic features of the reconstructions. Figure~\ref{subfig:Ini} shows the speed of sound reconstruction of the data simulated at \SI{0.5}{\milli\meter} resolution when reconstructed using a grid size of \SI{2}{\milli\meter}.  This reconstruction was computed using SLBFGS with source encoding gradient estimator $g_{\cE}$ for $128$ iterations (computing time~$\sim$\SI{30}{\minute}) and could be the multi-grid based initialization of a FWI run on the finer resolution. For this reason, we will call it $c_0^{ini}$ and take it as a reference point in the following studies. Note that here and in the following, we only reconstruct the speed of sound inside the breast region $\Gamma$, while the speed of sound outside is known and fixed. Figure~\ref{subfig:IniDiff} plots the reconstruction error $c_0^{ini} - c_0^{\dagger}$ locally, where $c_0^{\dagger}$ is the ground-truth speed of sound inside the breast. We will measure the global error of a variable $a$ by their relative $\ell_2$ distances to a reference $b$ computed as 
$\relleetwo(a, b) = \norm{a - b}_2 / \norm{b}_2$, and expressed in percent, e.g., we have that $\relleetwo(c_0^{ini}, c_0^\dagger) = 2.18\%$.

\begin{figure}[htb!]
   \centering
\subfloat[][\label{subfig:PreStudy32} 32 eval]{\includegraphics[height=0.4\textwidth]{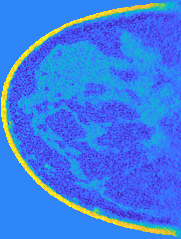}}
\hfill
\subfloat[][\label{subfig:PreStudy64} 64 eval]{\includegraphics[height=0.4\textwidth]{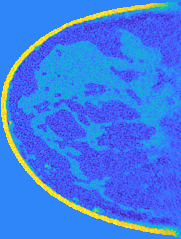}}
\hfill
\subfloat[][\label{subfig:PreStudy128} 128 eval]{\includegraphics[height=0.4\textwidth]{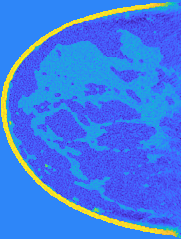}\hskip 2pt \includegraphics[height=0.4\textwidth]{ParulaNew1435-1665.pdf}}
\\
\subfloat[][\label{subfig:PreStudyDiff32} 32 eval]{\fbox{\includegraphics[height=0.4\textwidth]{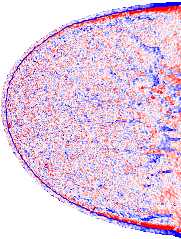}}}
\hfill
\subfloat[][\label{subfig:PreStudyDiff64} 64 eval]{\fbox{\includegraphics[height=0.4\textwidth]{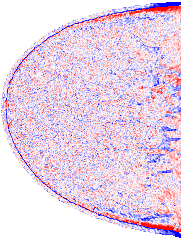}}}
\hfill
\subfloat[][\label{subfig:PreStudyDiff128} 128 eval]{\fbox{\includegraphics[height=0.4\textwidth]{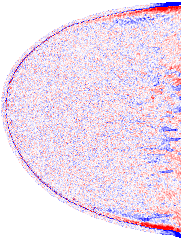}}\hskip 2pt \includegraphics[height=0.4\textwidth]{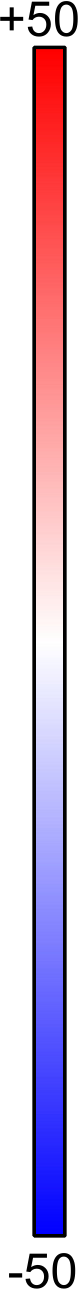}}
\caption{\protect\subref{subfig:PreStudy32}-\protect\subref{subfig:PreStudy128} $c_0^{rec}$ [m/s] computed by SLBFGS initialized at $c_0^{ini}$ with the source encoding gradient estimator $g_{\cE}$. \protect\subref{subfig:PreStudyDiff32}-\protect\subref{subfig:PreStudyDiff128} Corresponding error plots [m/s].}
   \label{fig:PreStudy1}
\end{figure}

\begin{figure}[htb!]
   \centering
\hfill
\subfloat[][\label{subfig:PreStudyDepth}]{\includegraphics[height=0.4\textwidth]{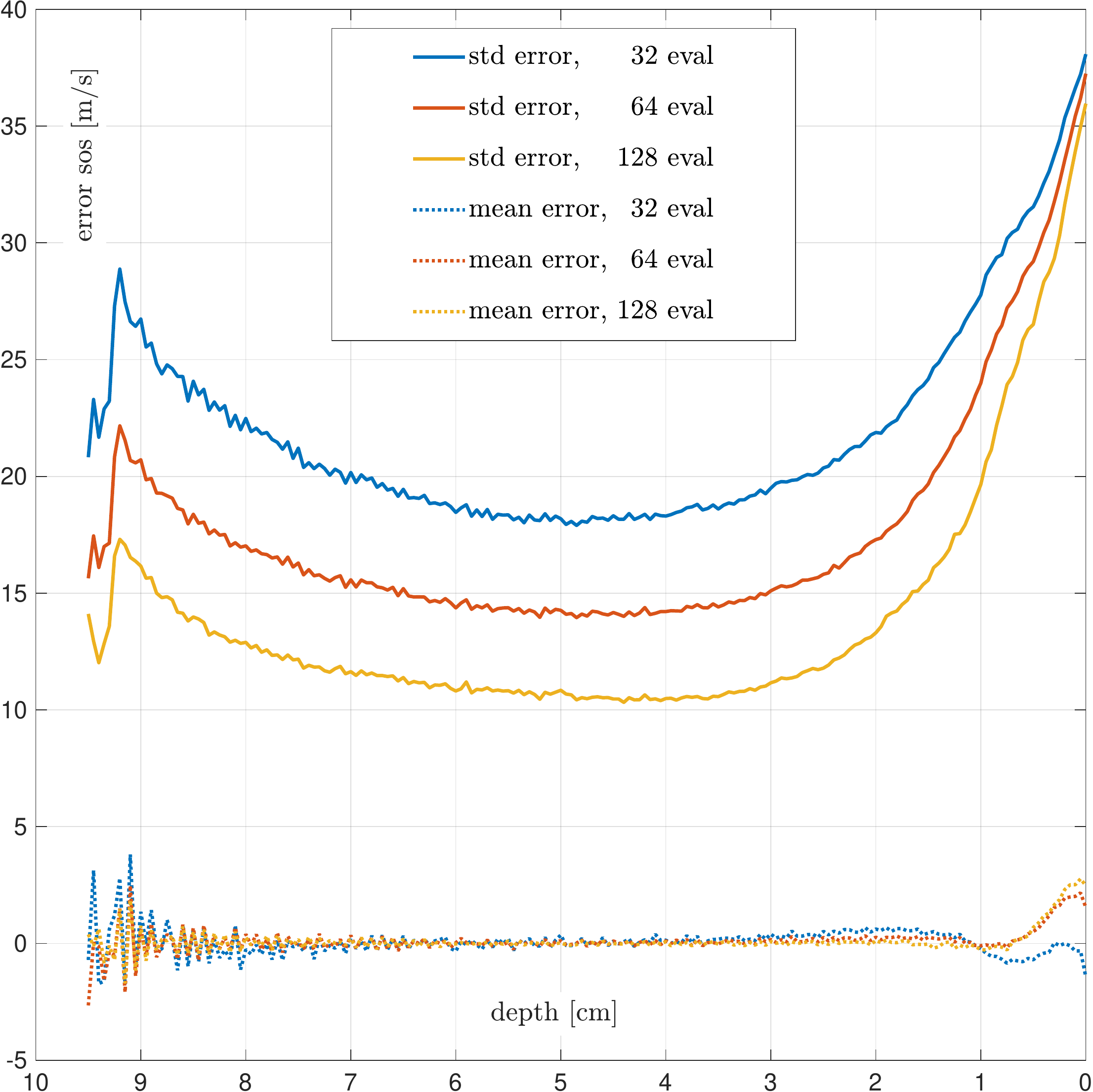}}
\hfill
\subfloat[][\label{subfig:TRGradDiff1tofull}]{\fbox{\includegraphics[height=0.4\textwidth]{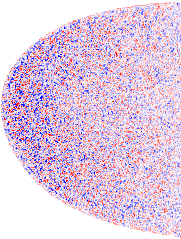}}\hskip 2pt \includegraphics[height=0.4\textwidth]{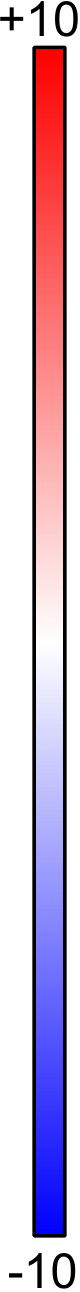}}
\hfill
\caption{\protect\subref{subfig:PreStudyDepth} Mean and standard deviation (std) of the reconstruction error $c_0^{rec} - c_0^\dagger$ vs. depth (= distance from chest wall) for the reconstructions shown in Figure~\ref{fig:PreStudy1}. \protect\subref{subfig:TRGradDiff1tofull} Difference between reconstruction using full gradient computation and TR based gradient with voxel layer of size 1.}
   \label{fig:PreStudy2}
\end{figure}

First, we run SLBFGS initialized at $c_0^{ini}$ with the source encoding gradient estimator $g_{\cE}$. Figure~\ref{fig:PreStudy1} shows the iterates after 32/64/128 gradient evaluations (\emph{eval}), the corresponding $\relleetwo$ errors are $1.53\%$/$1.29\%$/$1.07\%$. Note that we deliberately chose to not add measurement noise to the data (see above), the noisy appearance of the reconstructions are due errors of the gradient estimator (\emph{gradient noise}), a purely numerical phenomena that we need to highlight. On the difference plots in Figure~\ref{subfig:PreStudyDiff32}-\ref{subfig:PreStudyDiff128}, one can see that the error is larger around tissue interfaces and close to the chest wall. The former comes from the finite spatial grid spacing: the \textit{k}-space pseudospectral method cannot propagate the highest spatial frequencies present for a given spatial grid which means that we can at best hope to reconstruct a slightly smoothed version of $c_0^{\dagger}$. The higher errors close to the chest wall are caused by the insufficient sensor coverage of the half-spherical sensor array (cf. Figure~\ref{fig:Setup}). To illustrate this effect in more detail, we computed error statistics for each depth slice. For slices deeper than 2cm into the bowl, the errors are almost normally distributed, with zero mean and a standard deviation that decays as the SLBFGS advances as can be seen in Figure~\ref{subfig:PreStudyDepth}. We gain an accuracy of $\sim$\SI{5}{\meter / \second} with each doubling of number of iterations, i.e., the computational effort.

\subsection{Accuracy of Time Reversal Based Gradient Computation} \label{subsec:ValTRGrad}

First, we validate the TR based gradient computation described in Section \ref{subsec:TRGradient} when implemented using a \textit{k}-space pseudospectral method (cf. Appendix \ref{sec:kSpaceImp}). In principle, this schemes relies on regular spatial grids and is not well suited for representing arbitrarily shaped boundaries $\partial\Gamma$, cf. \eqref{eq:TRBwdWave}. To record or impose pressure values on  $\partial\Gamma$, we need to instead record and impose them on a layer of grid points lying on the inner surface of $\Gamma$. Determining the thickness of this layer is a trade-off: A thin layer leads to a small memory footprint when storing $p(x,t)$ on them, which was our original motivation to introduce the TR-based gradient computation. A large layer will lead to a more accurate reproduction of the forward field, i.e., $p^{\triangleleft}(x,t)$ is a better approximation of $p(x,T^{\triangleleft} - t)$. To examine the effects of this, we first compute the gradient \eqref{eq:Gradc0} at $c_0^{ini}$ for a single source storing the complete field $p(x,t)$ in $\Gamma \times [0,T]$. Then, we compute it using the TR formulation using boundary layers of thickness 1/2/4/8 voxels. Table \ref{table:TRGrad} shows the average relative errors of the TR gradients computed over 16 randomly picked sources and the working memory requirements in our setup. The errors decay very fast as a function of boundary layer thickness. However, we need to check that the errors in the TR gradients do not accumulate during an iterative reconstruction: Table \ref{table:TRGrad} shows the reconstruction errors when running SLBFGS with 32 gradient evaluations initialized at $c_0^{ini}$. Interestingly, even using a boundary layer of just one voxel, the error between reconstructions with normal gradient and TR-based gradient computations are small. Furthermore, Figure~\ref{subfig:TRGradDiff1tofull} shows that the error appears noise-like with no clear pattern visible. This suggests that the errors rather cancel than accumulate over the course of the iteration. The difference in computational time will be examined in Section \ref{subsec:GPUsCmp}. For the rest of the studies, a TR based gradient computation will be used with a boundary layer of 8 voxels. 

\begin{table}[tb!]
\centering
\begin{tabular}{|l|r|r|r|r|r|} 
 \hline
 & full & TR 1 & TR 2 & TR 4 & TR 8   \\ 
\hline
$\relleetwo$ gradients at $c_0^{ini}$ & 0 & $22.8085\%$  &  $5.6131\%$ & $0.8953\%$ & $0.1376\%$\\
memory requirements & \SI{98.00}{\giga\byte} & \SI{2.57}{\giga\byte} & \SI{5.95}{\giga\byte} & \SI{11.34}{\giga\byte} & \SI{21.75}{\giga\byte}\\
$\relleetwo(c_0^{rec}, c_0^{\dagger})$ & $1.5328\%$ & $1.5492\%$  &  $1.5401\%$ & $1.5319\%$ & $1.5328\%$\\
$\relleetwo(c_0^{rec, TR}, c_0^{rec, full})$ & $0\%$ & $0.2389\%$  &  $0.0681\%$ & $0.0178\%$ & $0.0038\%$\\
 \hline
\end{tabular}
\caption{Comparison of gradient computation storing the full forward field (\textit{full}) with TR based gradient computation using different number of boundary layer voxels.}
\label{table:TRGrad}
\end{table}

\subsection{Accuracy of Stochastic Gradient Optimization} \label{subsec:ValStochOpt}

\begin{figure}[tb!]
   \centering
\hfill
\subfloat[][\label{subfig:StochOptCmpEnergy}]{\includegraphics[height=0.4\textwidth]{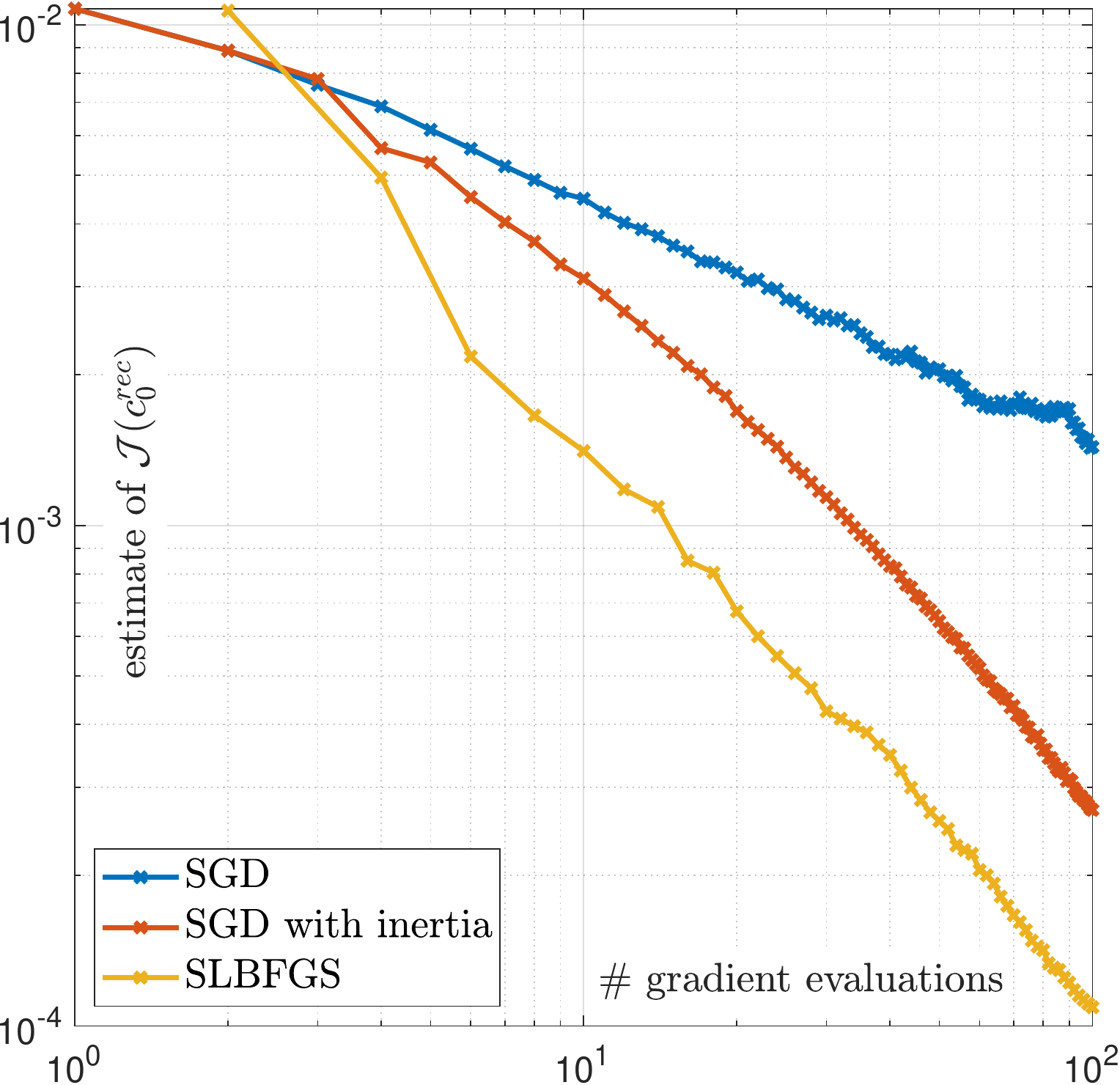}}
\hfill
\subfloat[][\label{subfig:StochOptCmpL2Err}]{\includegraphics[height=0.4\textwidth]{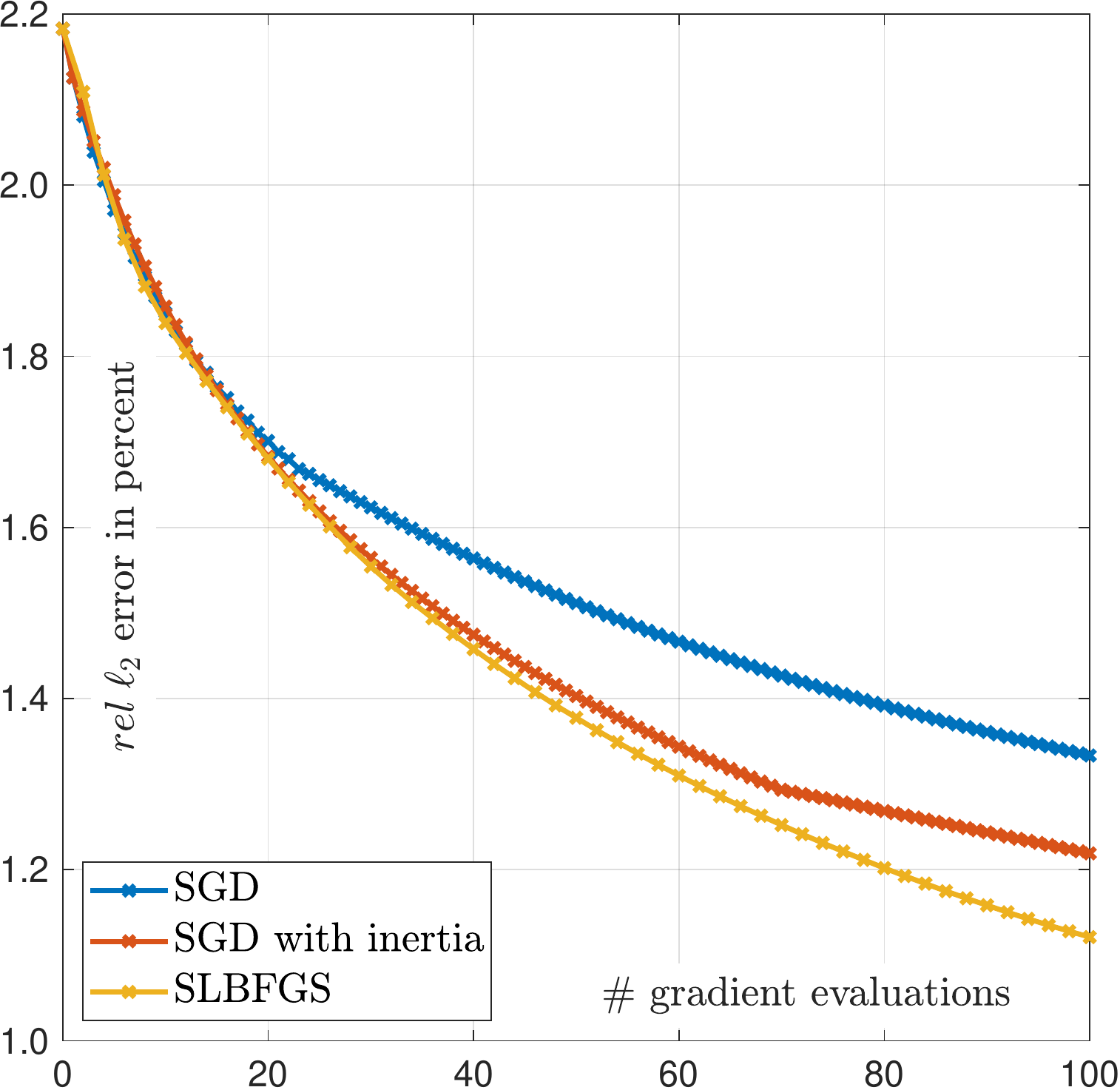}}
\hfill
\caption{Comparison of different stochastic gradient methods. \protect\subref{subfig:StochOptCmpEnergy} Estimated energy $\cJ(c_0^{rec})$ is plotted vs. the number of gradient evaluations. \protect\subref{subfig:StochOptCmpL2Err} $\relleetwo(c_0^{rec}, c_0^\dagger)$ is plotted vs. the number of gradient evaluations.}
   \label{fig:StochOptCmp}
\end{figure}

First, we demonstrate the impact the choice of the stochastic gradient method has and compare SLBFGS to plain SGD and SGD with inertia/momentum \cite{BoCuNo18} in Figure~\ref{fig:StochOptCmp}. All methods use source encoding as a gradient estimator and a maximum of $100$ gradient evaluations. Furthermore, once an increase in the estimated energy is detected, they do not return the iterates $u^k$ directly but a weighted average $\bar{u}^k := (\sum_l^k l^3 u^l) / \sum_l^k l^3$ to reduce variance, cf. Appendix \ref{sec:SLBFGS}. SLBFGS reaches the final $\relleetwo$ error of plain SGD already after $55$ evaluations, so twice as fast. This number increases to $76$ when we add inertia of $\beta = 0.5$ \cite{BoCuNo18} to SGD. However, we see that while SGD with inertia converges as fast as SLBFGS in the beginning, the decay levels off after $\sim60$ evaluations. 

\begin{figure}[tb!]
   \centering
\hfill
\subfloat[][\label{subfig:GradIni}]{\fbox{\includegraphics[height=0.4\textwidth]{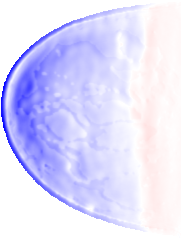}}\hskip 2pt \includegraphics[height=0.3\textwidth]{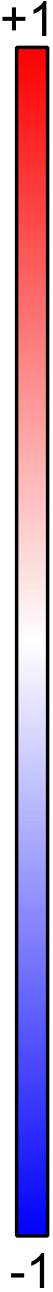}}
\hfill
\subfloat[][\label{subfig:StochGradEst1}]{\includegraphics[height=0.4\textwidth]{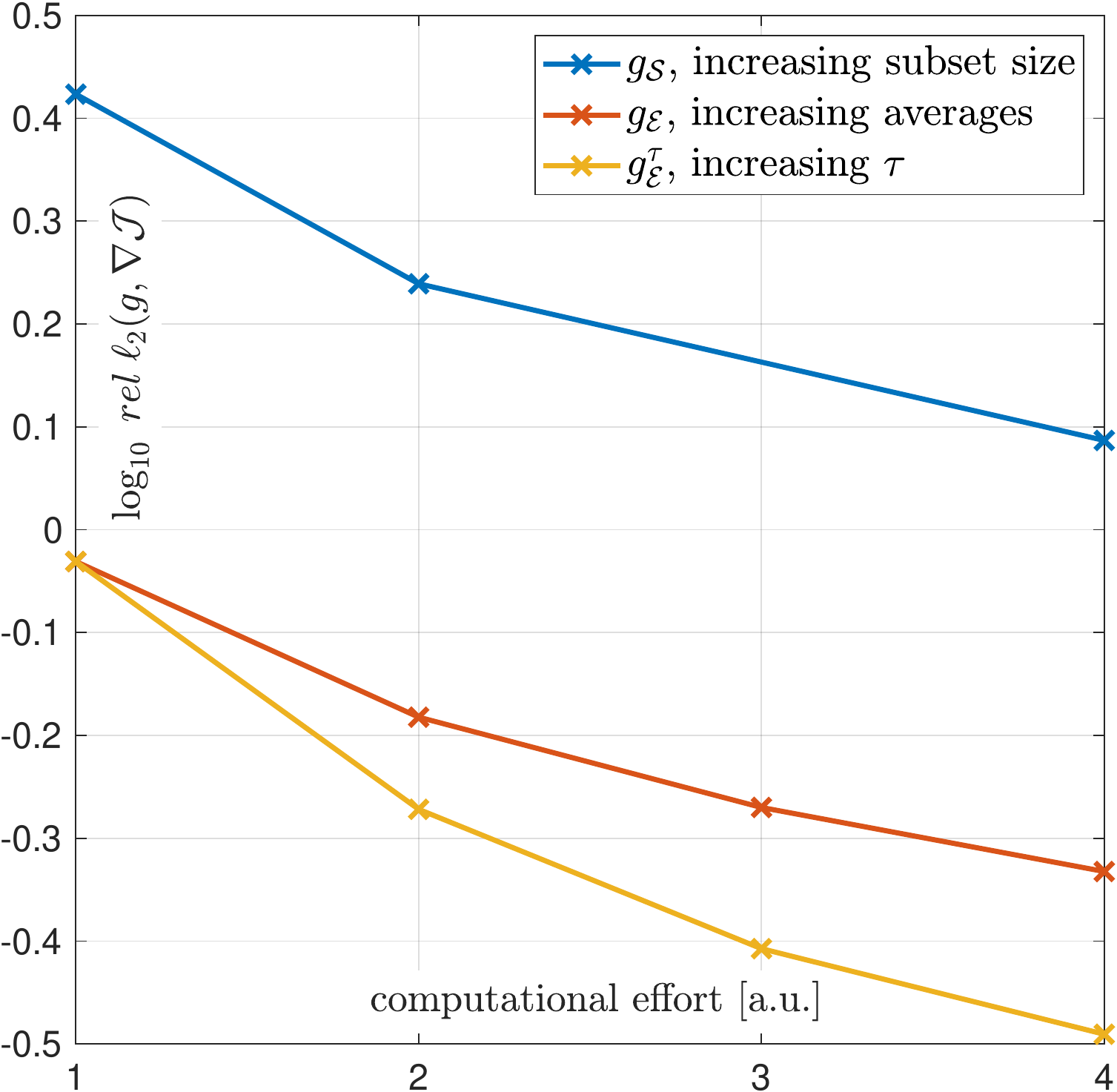}}
\hfill
\caption{\protect\subref{subfig:GradIni} Gradient $\nabla \cJ(c_0^{ini})$, normalized to $[-1,1]$. \protect\subref{subfig:StochGradEst1} Error of different gradient estimators and different ways to increase their precision by investing more computational effort, averaged over $16$ independent realizations. A computational effort of 1 corresponds to computing the gradient for a single source.}
   \label{fig:StochGradEst}
\end{figure}

\begin{figure}[tb!]
   \centering
\hfill
\subfloat[][\label{subfig:DelayedSELowAcc}]{\includegraphics[height=0.4\textwidth]{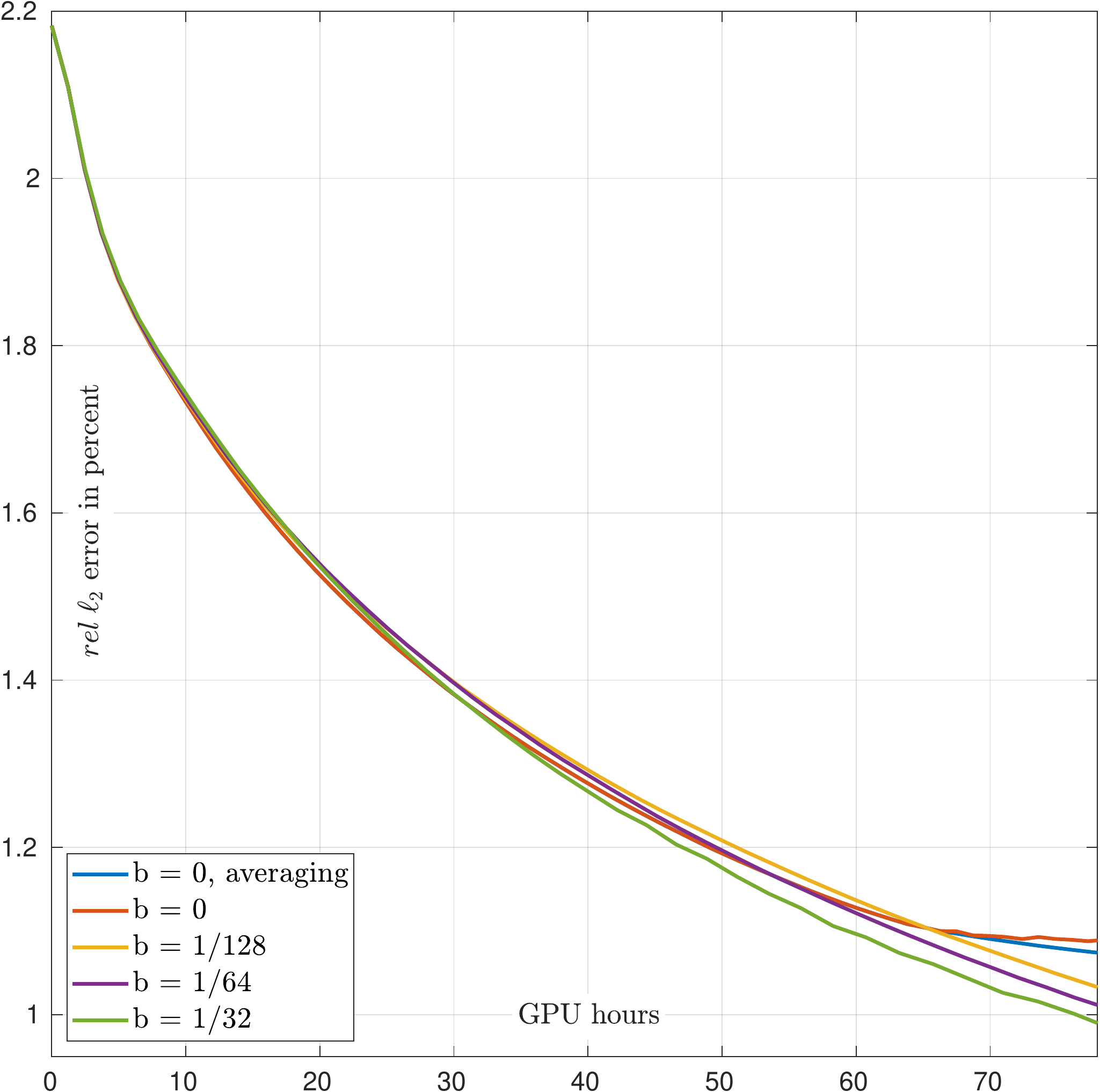}}
\hfill
\subfloat[][\label{subfig:DelayedSEHighAcc}]{\includegraphics[height=0.4\textwidth]{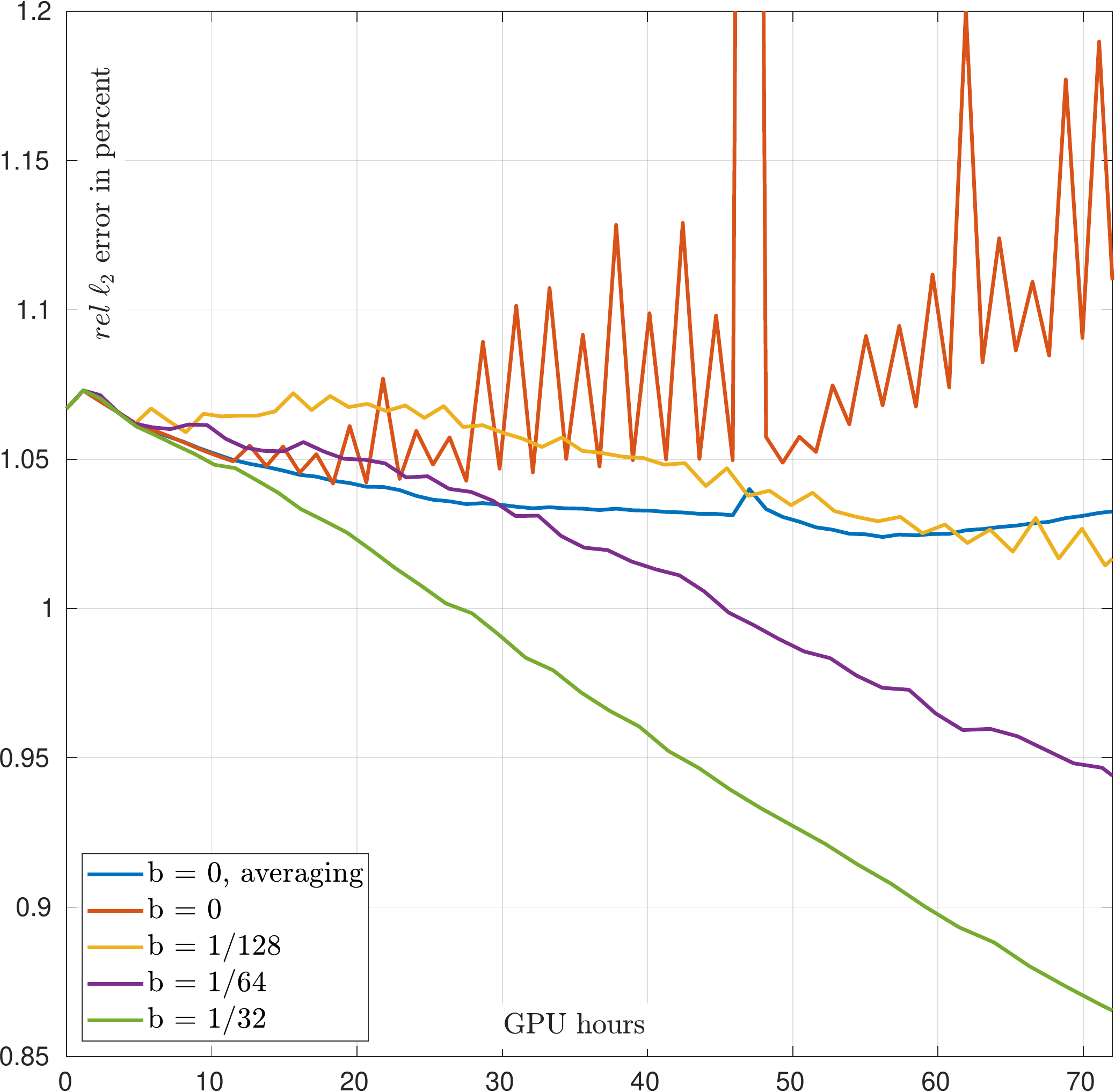}}
\hfill
\caption{Comparison of different delayed source encoding strategies. $\relleetwo(c_0^{rec}, c_0^\dagger)$ is plotted vs. computing time on a TITAN RTX GPU. All schemes are initialized with \protect\subref{subfig:DelayedSELowAcc} $c_0^{ini}$, cf. Figure \ref{subfig:Ini}; \protect\subref{subfig:DelayedSEHighAcc} the solution shown in Figure \ref{subfig:PreStudy128}.}
   \label{fig:DelayedSE}
\end{figure}

Second, we examine the impact of using different gradient estimators: Figure \ref{subfig:GradIni} shows the full gradient computed at $c_0^{ini}$, which took a week of computing time on a server with 4 GPUs (cf. Section \ref{subsec:GPUsCmp}). Then, we computed the different stochastic estimates of this gradient, $g_\cS$, $g_\cE$, $g_{\cE}^\tau$ (cf. Section \ref{subsec:SGDandSE}) and examined how their precision increases by investing more computational effort. Figure~\ref{subfig:StochGradEst1} shows the results averaged over $16$ i.i.d. gradient realizations. One can see that using source encoding ($g_\cE$,  $g_{\cE}^\tau$) indeed leads to a more accurate gradient estimation at the same computational effort. We will use source encoding for the rest of the studies. Also, we can see that using the novel delayed source encoding scheme $g_{\cE}^\tau$ not only gives a finer variance control compared to increasing the precision of $g_{\cE}$ via averaging, it also provides a more accurate gradient estimate with the same computational effort. Within a stochastic gradient method, the question is whether to do more steps with a less accurate but computationally cheaper gradient estimator vs. less steps with a more accurate but computationally more costly estimator. The former is typically more beneficial in the beginning of the iteration, while the latter becomes more beneficial towards the end of it \cite{BoCuNo18}. To illustrate this, we first initialized SLBFGS with $c_0^{ini}$ and then with the result shown in Figure \ref{subfig:PreStudy128}, which corresponds to 128 gradient evaluations (roughly 4 days computing time on a TITAN RTX GPU). In both cases, we ran it for another 72 hours using $g_{\cE}^\tau$ with $\tau = b\, (i - 1)$ for different values of $b$ and with and without starting iterate averaging after an increase in the energy estimate in the case $b = 0$ (so no delay). Figure \ref{subfig:DelayedSELowAcc} shows that in the low accuracy regime, all methods perform rather similar for the first 24 hours. In the case $b=0$, iterate averaging is only activated after 66 hours (the point where the first two plots diverge). The delayed source encoding methods $b > 0$ eventually perform better than $b=0$, but are slower at the start. Figure \ref{subfig:DelayedSEHighAcc} shows the results for the high accuracy regime: While using iterate averaging stabilizes the iteration for $b=0$, using increasing values of $b$ leads to faster convergence in terms of computing time.

\subsection{Multi-GPU Acceleration} \label{subsec:MultipleGPUs}

\begin{figure}[tb!]
   \centering
\hfill
\subfloat[][\label{subfig:MultGPUErr}]{\includegraphics[height=0.4\textwidth]{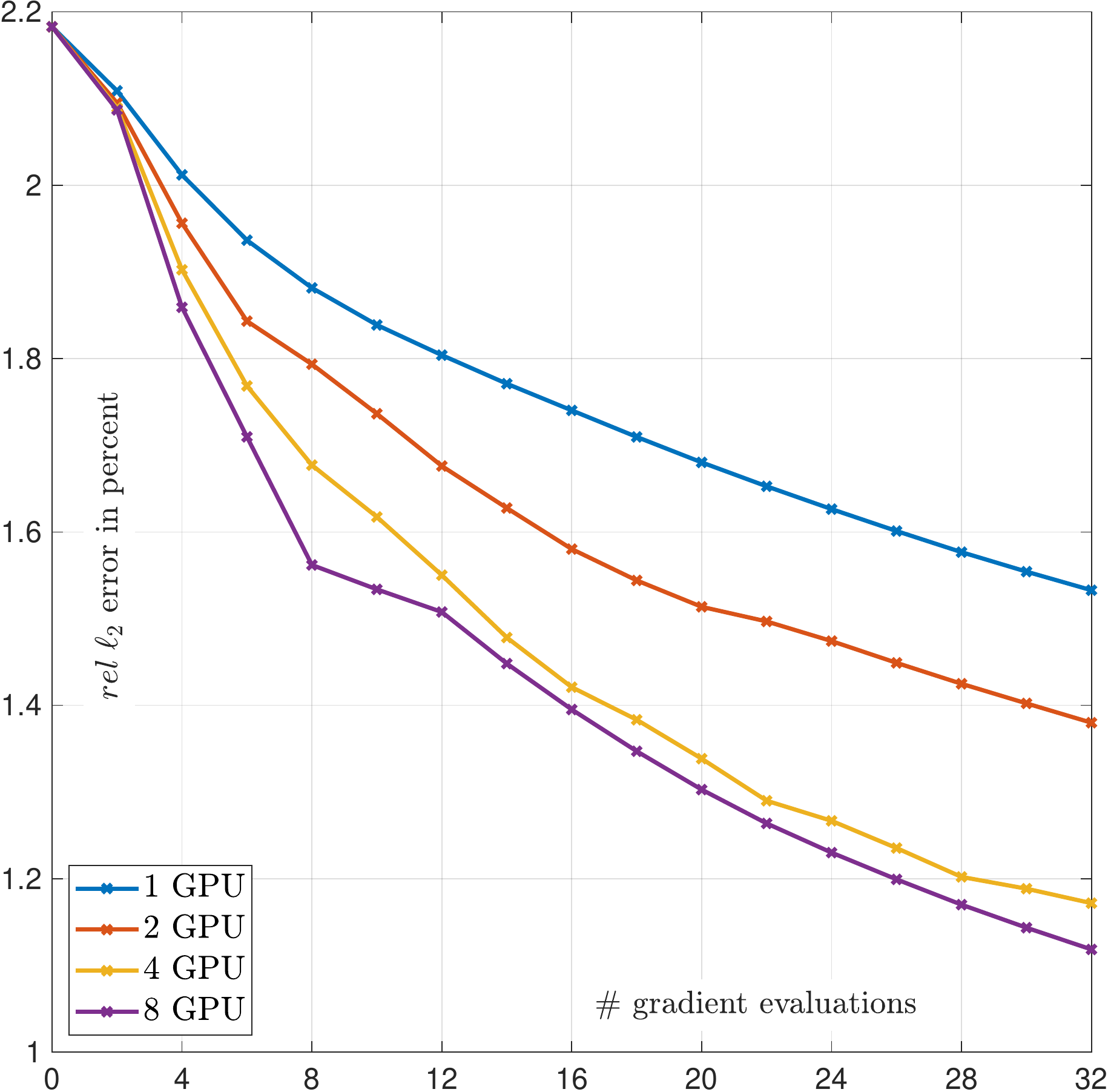}}
\hfill
\subfloat[][\label{subfig:MultGPUDepth}]{\includegraphics[height=0.4\textwidth]{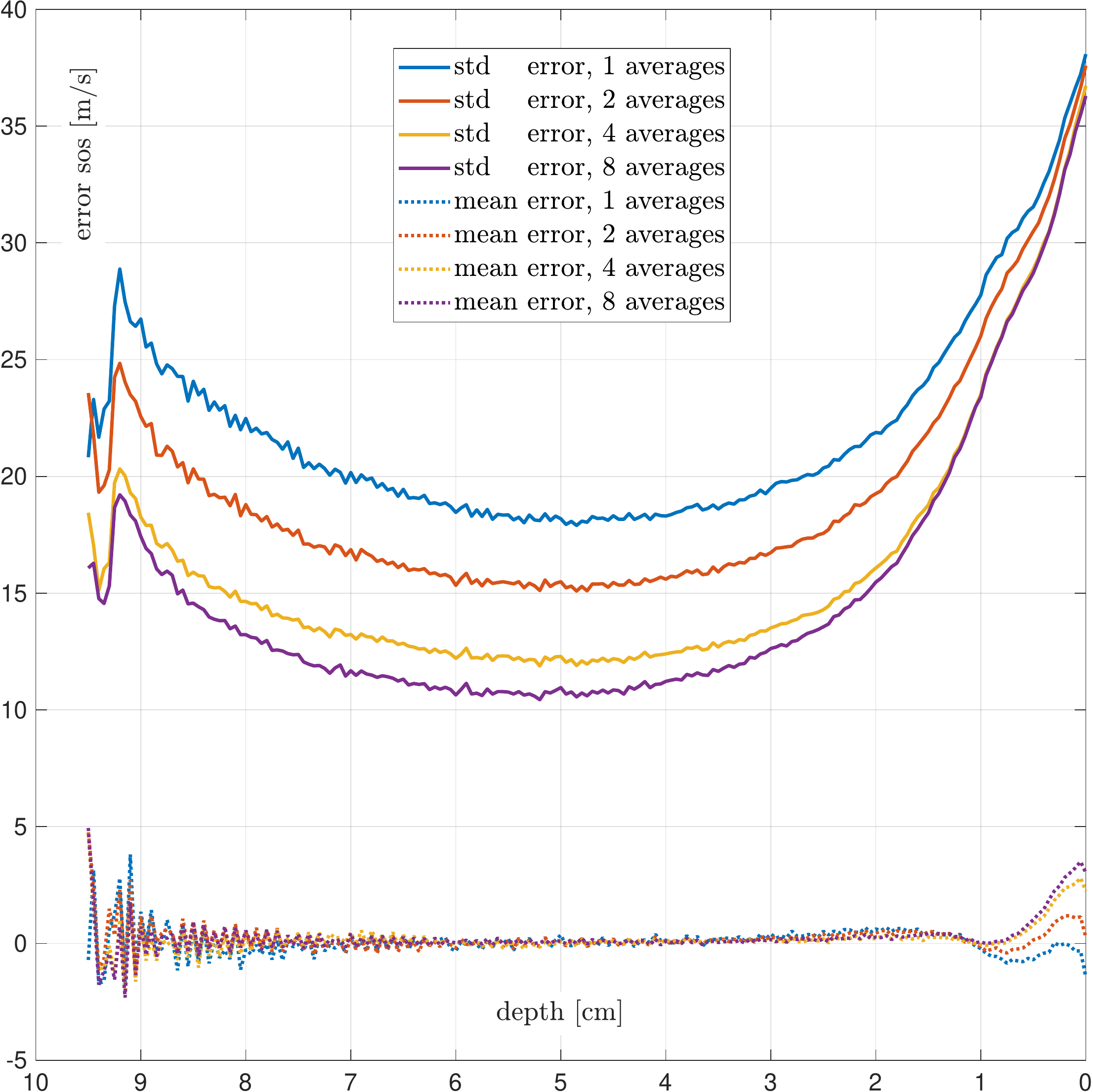}}
\hfill
\caption{Effect of utilizing multiple GPUs to decrease the variance of the SE gradient estimator by averaging independent realizations. \protect\subref{subfig:MultGPUErr}  $\relleetwo(c_0^{rec}, c_0^\dagger$) is plotted vs. the number of gradient evaluations. \protect\subref{subfig:MultGPUDepth} Mean and standard deviation (std) of the reconstruction error $c_0^{rec} - c_0^\dagger$ vs. depth, cf. Figure~\ref{subfig:PreStudyDepth}.}
   \label{fig:MultGPU}
\end{figure}

Computing platforms with multiple GPUs can be used in different ways to accelerate the image reconstruction. In the TR based gradient computation described in Section \ref{subsec:TRGradient}, the adjoint and TR wave simulation can be run on two separate GPUs in parallel, which would lower the computational cost to that of the standard gradient computation again. To reconstruct 3D images on higher resolutions, one can implement 3D domain decomposition methods to distribute the computations over several GPUs \cite{NaJaReTr12,TrVaJa16,TrVaJa18}. A straightforward way which does not need sophisticated implementation is to average statistically independent gradient estimates each computed on a different GPU. This reduces the variance of the gradient estimator at least by a factor corresponding to the number of GPUs, cf. Section \ref{subsec:SGDandSE}: Figure~\ref{fig:MultGPU} shows the reconstruction errors of SLBFGS initialized at $c_0^{ini}$ when using a growing number of GPUs in this way (using source encoding without time shifting). One can see that the convergence is in general faster. However, Figure~\ref{subfig:MultGPUErr} highlights that due to the non-convexity of $\cJ(u)$, there is not a simple relationship between the accuracy of the gradient estimator and the convergence curve. The reconstruction error obtained after running SLBFGS using 2/4/8 GPUs for 32 gradient evaluations is the same as running it on a single GPU for 51/89/103 gradient evaluations. Comparing Figures~\ref{subfig:PreStudyDepth} and ~\ref{subfig:MultGPUDepth} shows also the decrease in standard deviation saturates with a growing number of GPUs.   

\subsection{Comparison of Different GPU Architectures} \label{subsec:GPUsCmp}

Table \ref{table:GPU} compares the performance of different GPUs. One can see that the run times of normal and TR-based gradient computations are more alike than expected based on the number of computational operations (\textit{flops}). One reason is that the huge memory footprint of the normal gradient computation requires a lot of working memory operations, which are no longer negligible compared to the computational operations. 

\begin{table}[tb!]
\centering
\begin{tabular}{|l|r|r|r|r|} 
 \hline
 & Tesla P40 & TITAN RTX & GeForce RTX 2080 Ti & GeForce GTX 1070    \\ 
\hline
memory & \SI{24}{\giga\byte} & \SI{24}{\giga\byte} & \SI{11}{\giga\byte} & \SI{8}{\giga\byte}    \\  
fwd     & 19m 54s  & 10m 59s & 11m 42s & 30m 27s   \\ 
grad    & 46m 9s   & 26m 44s & 27m 34s & n/a \\ 
TR grad & 61m 48s  & 34m 22s & 36m 15s & 1h 34m \\ 
 \hline
\end{tabular}
\caption{Comparison of different Nvidia GPUs. Shown are the computation times for simulating data generated by a single source (\textit{fwd}), computing a gradient for it via the standard approach \eqref{eq:Gradc0} (\textit{grad}) and the TR based computation (cf. Section \ref{subsec:TRGradient}) (\textit{TR grad}). N/A indicates that the computer the GPU was installed in did not have insufficient working memory (RAM).}
\label{table:GPU}
\end{table}

\subsection{Acceleration by Coarse-To-Fine Multigrid Initilization} \label{subsec:ValMultGrid}

\begin{figure}[tb!]
   \centering
\hfill
\subfloat[][\label{subfig:MultGridErr}]{\includegraphics[height=0.4\textwidth]{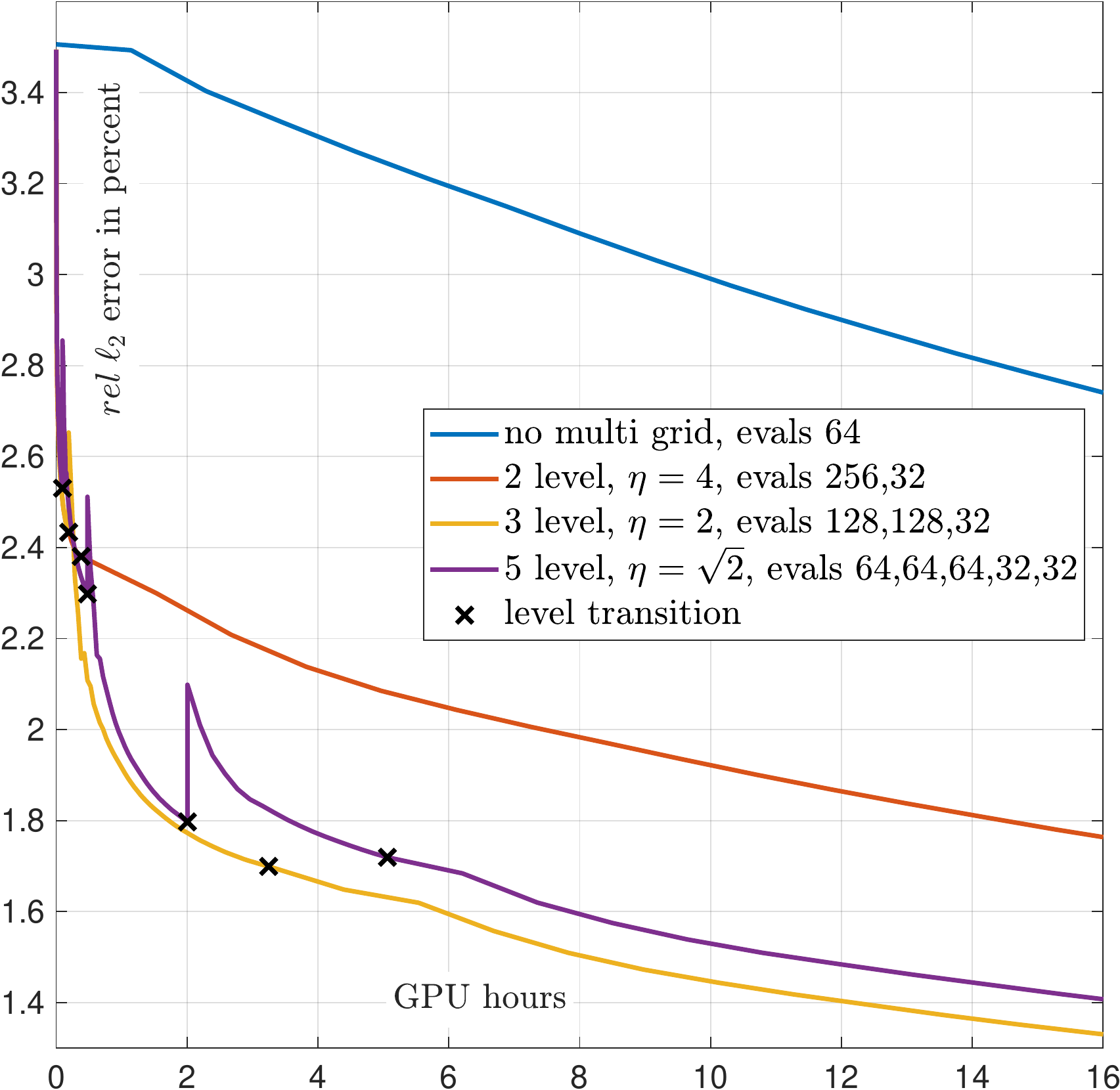}}
\hfill
\subfloat[][\label{subfig:PreconErr}]{\includegraphics[height=0.4\textwidth]{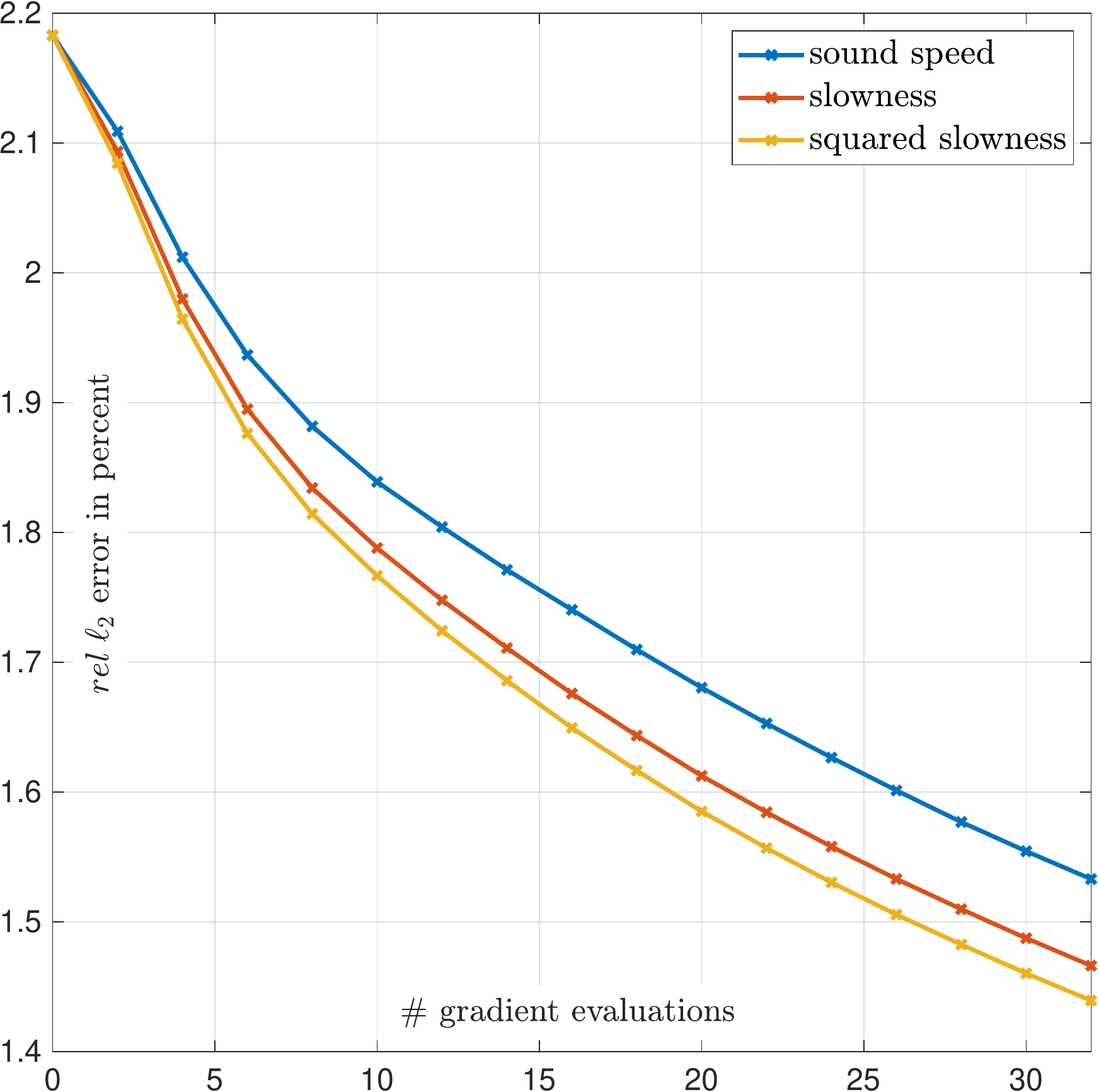}}
\hfill
\caption{\protect\subref{subfig:MultGridErr} Comparison of using multi-grid initilizations: All schemes are initialized with a uniform background ($c_0 = $\SI{1500}{\meter / \second}). A different number of gradient evaluations (evals) are run on each level, the TR-based approach is only used on the finest level ($dx =$ \SI{0.5}{\milli \meter}). In principle, the vertical axis displays flops which scale as $\eta^4$ (cf. Section \ref{subsec:MultGrid}) where $\eta$ is the coarsening factor between levels. For a better interpretability, we re-scaled flops to GPU hours by using the TR-gradient computation time of the TITAN RTX on the finest level (cf. Table \ref{table:GPU}). The black crosses mark where a scheme transitions between levels. The $\relleetwo$ error of solutions on coarser grids than $dx =$ \SI{0.5}{\milli \meter} is measured by interpolating them first. Note that each level transition involves a smoothing which can lead to an increase in  $\relleetwo$ error. \protect\subref{subfig:PreconErr} $\relleetwo(c_0^{rec}, c_0^\dagger$) vs. the number of gradient evaluations for SLBFGS initialized at $c_0^{ini}$ and different preconditioning techniques. }
   \label{fig:MultGridPrecon}
\end{figure}

Next, we illustrate that a coarse-to-fine initialization strategy can speed up the convergence and avoid local minima as described in Section \ref{subsec:MultGrid}. We use SLBFGS with source encoding without time shifting and initialize it with the sound speed value of water (\SI{1500}{\meter / \second}) everywhere. The first scheme runs on the finest grid ($dx = 0.5$\si{\milli\meter}) only. The second scheme has two levels: It starts on $dx = 2$\si{\milli\meter} before switching to $dx = 0.5$\si{\milli\meter}, so $\eta = 4$. The third scheme has three levels $dx = 2, 1, 0.5$\si{\milli\meter} ($\eta = 2$), and the forth $dx = 2, 1.41, 1, 0.70, 0.5$\si{\milli\meter} ($\eta = \sqrt{2}$). Figure~\ref{subfig:MultGridErr} shows the error plots vs computational effort and shows how essential it is to use multi-grid schemes in UST: When the 2-level scheme switches to the finest level after 23 minutes, it has already reached a lower error than the single-level scheme will reach after 16 hours. 

\subsection{Acceleration by Preconditioning} \label{subsec:ValPrecon}

Figure \ref{subfig:PreconErr} shows that using the non-linear slowness and squared-slowness preconditioning  transformations discussed in Section \ref{subsec:MultGrid} do indeed speed up the convergence and provide clear advantages over the standard parameterization. All spatial error plots clearly indicate that the errors increase with decreasing depth, i.e., closer to the chest wall where the sensor coverage is getting worse, cf. Figures \ref{fig:Setup}-\ref{fig:PreStudy2}. This observation suggests that introducing a local weighting $w(x)$ that depends on depth or sensor coverage could be beneficial. We tried different variants of such weightings but could not find significant differences. One explanation for this is that quasi-Newton schemes like SLBFGS have low sensitivity to linear preconditioning transformations.

\section{Full Inversion Scheme Demonstration}  \label{sec:FullDemo}

In this section, we combine all the techniques described in the previous section to demonstrate the results one can obtain within a computing time limit of 24 hours. While this limit may seem arbitrary at this point, it is motivated by the perspective of adding UST to clinical work-flows for breast-cancer diagnosis and treatment planning: After a positive screening result, several examinations have to be scheduled and carried out and before a treatment decision can be made. While health care providers try to shorten the overall time for this as much as possible, it can rarely be done in a single day\cite{SiVoGrVeMaBrVeMuGeBr20}. Within this context, a computing time of 24 hours would allow UST to fit into clinical trajectories. Similar considerations hold for the computational resources needed: A main advantage of UST over MRI could be lower operating costs for a single examination. Having to use or rent a large computational cluster over several days to compute the results of a single examination may diminish this advantage.

\begin{figure}[b!]
   \centering
\hfill
\subfloat[][\label{subfig:RecStart}]{\includegraphics[height=0.4\textwidth]{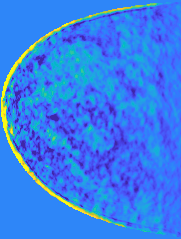}
\hskip 2pt \includegraphics[height=0.4\textwidth]{ParulaNew1435-1665.pdf}}
\hfill
\subfloat[][\label{subfig:RecDepth}]{\includegraphics[height=0.4\textwidth]{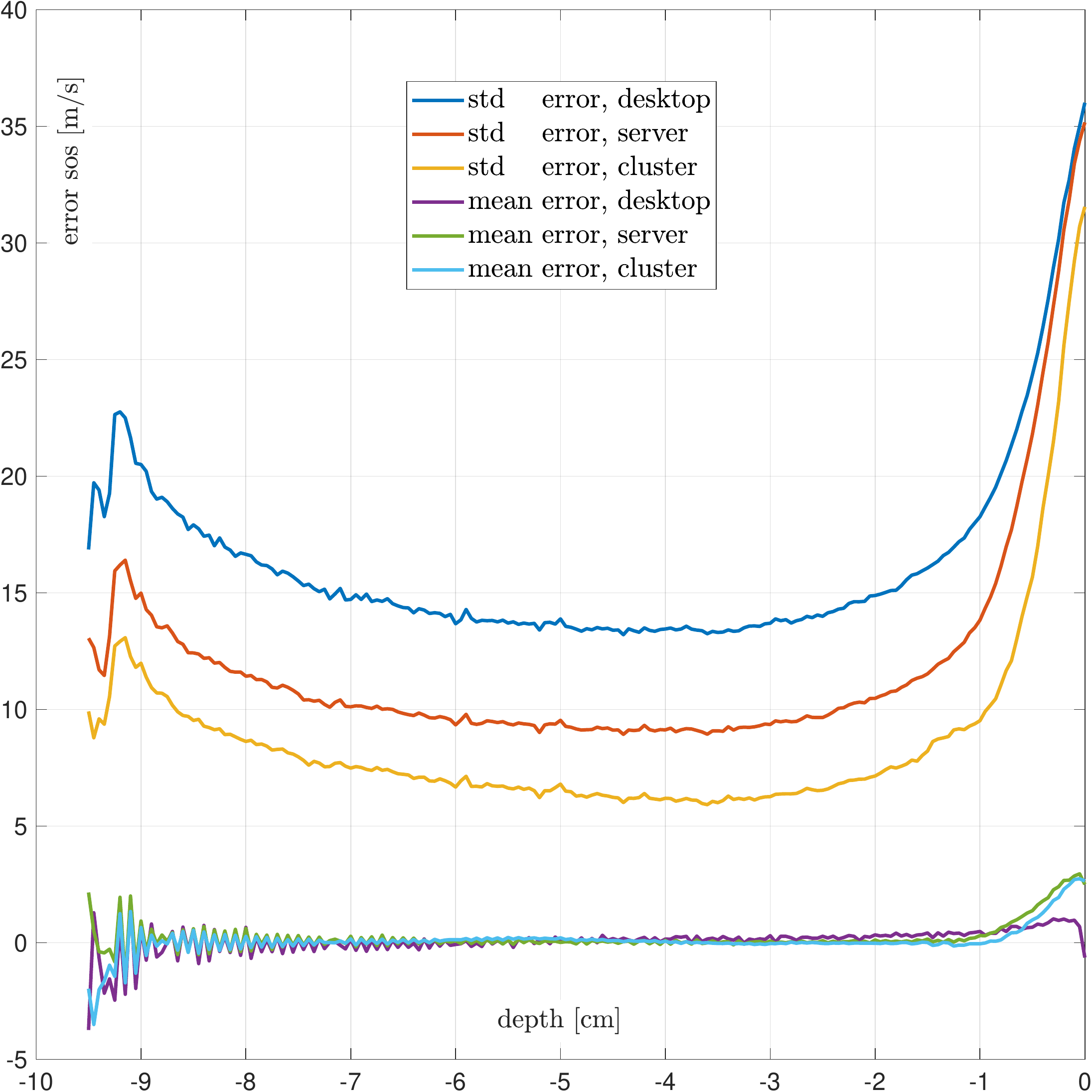}}
\hfill
\caption{\protect\subref{subfig:RecStart} Reconstruction a reconstruction method similar to the one used \cite{PeHeUdMiCoTr17} within \SI{24}{\hour} computing time on a single GPU. \protect\subref{subfig:RecDepth}  Mean and standard deviation (std) of the reconstruction error $c_0^{rec} - c_0^\dagger$ vs. depth, cf. Figure~\ref{subfig:PreStudyDepth}.}
   \label{fig:FullDemoDepth}
\end{figure}

\begin{figure}[tb!]
   \centering
\hfill
\subfloat[][\label{subfig:RecDesktop}]{\includegraphics[height=0.4\textwidth]{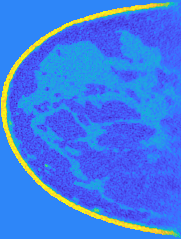}}
\hfill
\subfloat[][\label{subfig:RecServer}]{\includegraphics[height=0.4\textwidth]{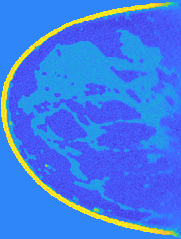}}
\hfill
\subfloat[][\label{subfig:RecCluster}]{\includegraphics[height=0.4\textwidth]{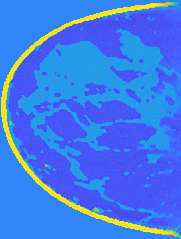}
\hskip 2pt \includegraphics[height=0.4\textwidth]{ParulaNew1435-1665.pdf}}
\hfill\\
\hfill
\subfloat[][\label{subfig:RecDesktopDiff}]{\fbox{\includegraphics[height=0.4\textwidth]{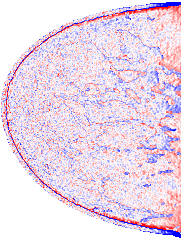}}}
\hfill
\subfloat[][\label{subfig:RecServerDiff}]{\fbox{\includegraphics[height=0.4\textwidth]{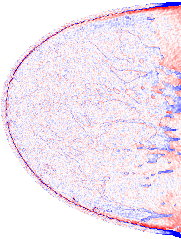}}}
\hfill
\subfloat[][\label{subfig:RecClusterDiff}]{\fbox{\includegraphics[height=0.4\textwidth]{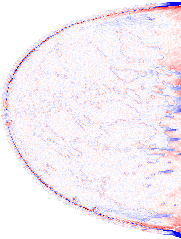}}\hskip 2pt \includegraphics[height=0.4\textwidth]{Red2BlueColorbar50.pdf}}
\hfill
\caption{Demonstration full inversion scheme. Shown are results computed within \SI{24}{\hour} using \protect\subref{subfig:RecDesktop} a desktop PC with one GeForce RTX 2080 Ti GPU,  \protect\subref{subfig:RecServer} a server with 4 GeForce RTX 2080 Ti GPUs, \protect\subref{subfig:RecCluster} a cluster with 16 TITAN RTX GPUs. \protect\subref{subfig:RecDesktopDiff}-\protect\subref{subfig:RecClusterDiff} corresponding error plots. While reconstruction \protect\subref{subfig:RecDesktop} still shows some stochastic gradient noise, all main anatomical features are visible and there is no systematic quantitative error.} 
   \label{fig:FullDemo}
\end{figure}

To appreciate the computational advantages described in this work, we first show the reconstruction result of applying a computational technique similar to the one used in \cite{PeHeUdMiCoTr17} in Figure~\ref{subfig:RecStart}. One can see that even after 24 hours the result can not be used to reliably identify anatomic structures. Then, we use
\begin{itemize}
    \item TR based gradient computation with a $8$ voxel boundary layer as examined in Section \ref{subsec:ValTRGrad}.
    \item SLBFGS with source encoding and weighted iterate averaging as examined in Section \ref{subsec:ValStochOpt}. 
    \item Multigrid initialization with $dx = 2, 1, 0.5$\si{\milli\meter} ($\eta = 2$) and $128/128/32$ gradient evaluations as examined in Section \ref{subsec:ValMultGrid}.
    \item Squared slowness preconditioning as examined in Section \ref{subsec:ValPrecon}.
\end{itemize}
The results of computing on a single GeForce RTX 2080 Ti installed in a conventional desktop PC (total price $\sim5$K\euro) are shown in Figures \ref{subfig:RecDesktop}/\ref{subfig:RecDesktopDiff}. While one can still perceive stochastic gradient noise and structured errors in the reconstruction, it already provides a useful high-resolution estimate of the underlying anatomy. Figures   
\ref{subfig:RecServer}/\ref{subfig:RecServerDiff} and \ref{subfig:RecCluster}/\ref{subfig:RecClusterDiff} show how these results improve using a server with 4 GeForce RTX 2080 Ti (total price $\sim20$K\euro) or a cluster with 16 TITAN RTX (total price $\sim100$K\euro). Figure~\ref{subfig:RecDepth} shows the corresponding depth distribution of the errors. At this point, we remind the reader that the images shown only illustrate a sub-set of the whole reconstruction domain, cf. Figure~\ref{fig:Setup}. 

As discussed earlier, we did not add measurement noise to the data in our studies to allow for an easier visual assessment of the convergence of the different computational schemes as the stochastic gradient error also manifests as noise-like image features. To round off the numerical studies, Appendix \ref{sec:Noise} examines the sensitivity to noise.

\clearpage

\section{Discussion and Outlook} \label{sec:DisOut}

Our work was motivated by the specifications of the PAMMOTH scanner that is currently being developed and tested \cite{PammothWebsite}. The main contribution was to combine a number of techniques to design a comprehensive computational strategy to realize high resolution 3D TD-FWI for UST of the human breast with moderate computational resources within a day of computing time. The results were validated with extensive numerical proof-of-concept studies on different computational platforms.

\paragraph{TR-based gradient computation} For the biggest breast cup designed for PAMMOTH, storing $p(x,t)$ throughout the breast would require  \SI{146.81}{\giga\byte} at \SI{0.5}{\milli\meter} spatial resolution. This memory footprint rules out keeping it in the memory of any currently available GPU. If kept in CPU memory, it severely limits the number of parallel source simulations one can perform, e.g., running parallel simulations on a server with 16 GPUs like showcased in Section \ref{sec:FullDemo} would require more than \SI{2}{\tera\byte} CPU memory. In addition, memory operations start to dominate the total computation time (cf. Table \ref{table:GPU}). \\ 
The time-reversal-based gradient computation presented in Section \ref{subsec:TRGradient} can overcome this bottleneck by shifting the storage from the breast volume to its boundary: With a moderate increase of computation time (cf. Table \ref{table:GPU}), the memory requirement can be lowered by a factor $28.79$ to as little as \SI{5.10}{\giga\byte} without compromising accuracy of the final solution by more than $0.25\%$, cf. Section \ref{subsec:ValTRGrad}. This allows us to utilize multiple GPUs without memory restrictions, which is examined in detail in Section \ref{subsec:MultipleGPUs} and demonstrated in Section \ref{sec:FullDemo}, cf. Figures \ref{fig:MultGPU}, \ref{fig:FullDemo}. 

\paragraph{Stochastic optimization} The second crucial ingredient of our strategy is to combine novel stochastic optimization techniques with efficient stochastic gradient estimators and preconditioning. In Section \ref{subsec:ValStochOpt}, we demonstrate the benefits of SLBFGS over plain SGD and show that source encoding schemes lead to much more accurate gradient estimators at the same computational cost as sub-sampling methods, cf. Figure \ref{fig:StochGradEst}. Furthermore, we presented a novel source encoding scheme that exploits the time-invariance of the wave equation by introducing random time delays to increase the accuracy of the gradient estimate at the expense of longer computation times. This offers a much finer variance control than estimate averaging as well as being more accurate, cf. Figure \ref{subfig:StochGradEst1}. Embedded in SLBFGS, we showed that using delayed source encoding will be beneficial over iterate averaging in the high accuracy regime, beyond the fixed computation time budget of $24$ hours we set ourselves in this work, cf. Figure \ref{fig:DelayedSE}. Future work will examine this topic more carefully, e.g., by developing adaptive delay choice rules and by examining how to balance delays to optimally utilize heterogeneous multi-GPU computing environments. The encoding scheme presented here relies on the transient nature of the wave fields to reduce cross-talk. It would be interesting to compare it to the approach presented in \cite{BaTr20}, which obtains impressive results by encoding a harmonic steady state instead. Finally, one can observe that the error of the source encoding gradient estimates appears to be noise-like (cf. Figures \ref{fig:PreStudy1}, \ref{fig:FullDemo}). As such, it is worth examining whether coupling source encoding with preconditioning that performs image denoising will further accelerate the convergence.   

\paragraph{Multigrid schemes} The results presented in Section \ref{subsec:ValMultGrid} reveal that using coarse-to-fine multi-grid schemes are another essential component of our strategy: Even the simple multi-grid initialization technique presented here accelerates the overall convergence massively. In the future, we will examine more sophisticated multi-grid schemes.

\paragraph{Further Extensions} 

\begin{itemize} 
\item The methods presented here were exclusively implemented in \textsc{Matlab}, using the GPU support offered by the Parallel Computing Toolbox. Implementing it using optimized CUDA code will lead to further performance gains \cite{NaJaReTr12,TrVaJa16,TrVaJa18}.
\item While we presented the material in terms of a general acoustic parameter $u$, we only showcased the methods for $u = c_0$. The extension to acoustic density $\rho_0$ and attenuation $(\alpha_0, y)$ is described in Appendix \ref{sec:DenAbsGrad} and will be examined numerically in future work.
\item In this numerical proof-of-concept study, we used idealized assumptions on the transducer locations and properties. In forthcoming work, we will address realistic transducer and scan modeling.  
\item Finally, we will demonstrate the application of this work to real data from experimental phantoms, healthy volunteers and breast-cancer patients acquired by the PAMMOTH scanner in the future. 
\end{itemize}

\section{Conclusions} \label{sec:Con}

In this work, we described and evaluated a comprehensive computational strategy for high resolution 3D TD-FWI for UST of the human breast. With moderate computational resources, accurate results can be obtained within less than a day of computation time. The extent of the improvements presented here can best be appreciated by comparing Figure \ref{subfig:RecStart}, which corresponds to our starting point \cite{PeHeUdMiCoTr17} with Figure \ref{fig:FullDemo}, which showcases our current results. Similarly promising results were presented in \cite{BaTr20}, which was developed in parallel to the work described here. While the computing time and costs of our scheme are still too high for screening applications, it may well already fit into clinical trajectories for diagnosis and therapy planing.

\subsection*{Acknowledgements} \label{sec:Ack}
This work was supported by the Netherlands Organisation for Scientific Research (NWO) project number 613.009.106 and by the European Union’s Horizon 2020 Research and Innovation program H2020 ICT 2016-2017 under Grant agreement No. 732411, which is an initiative of the Photonics Public Private Partnership.


\appendix

\section*{Appendix}

\section{Extension to Density and Absorption} \label{sec:DenAbsGrad}

This section extends the formulas in Section \ref{subsec:FWI} and lists the remaining derivatives. First, including absorption with $L \ne 0$, we get
\begin{equation}
\frac{\partial L}{\partial c_0} = - 2 \alpha_0 (y-1) c_0^y \frac{\partial}{\partial t} \left(- \nabla^2 \right)^{\tfrac{y}{2}-1} + 2 \alpha_0 y c_0^{y-1} \frac{\partial}{\partial t} \left( - \nabla^2  \right)^{\tfrac{y+1}{2}-1}
\end{equation} 
as additional terms for the derivative with respect to $c_0$. See \cite{PeHeUdMiCoTr17} for a discussion of their relative importance. For $u = \rho_0$ we get 
\begin{equation}
\frac{\partial A}{\partial \rho_0} = \left(\nabla \cdot \frac{1}{\rho_0} - (\rho_0 \nabla \cdot) \frac{1}{\rho_0^2}\right) \nabla
\end{equation}
For $u = \alpha_0$ we get 
\begin{equation}
\frac{\partial A}{\partial \alpha_0} = \frac{\partial L}{\partial \alpha_0} = - 2 c_0^{y-1} \frac{\partial}{\partial t} \left( - \nabla^2  \right)^{\tfrac{y}{2}-1} + 2 c_0^y \tan(\pi y/ 2) \left( - \nabla^2  \right)^{\tfrac{y+1}{2}-1}
\end{equation}
For $u = y$, we remark that functions of the self-adjoint operator $-\nabla^2$, are understood in terms of its spectral decomposition, i.e., in Fourier space (cf. \cite{TrCo11b}) and get
\begin{multline}
    \frac{\partial A}{\partial y} = \frac{\partial L}{\partial y} = 
    \tau \log(c_0) \frac{\partial}{\partial t} \left(- \nabla^2 \right)^{\tfrac{y}{2}-1} + \frac{\tau}{2} \frac{\partial}{\partial t} \left(- \nabla^2 \right)^{\tfrac{y}{2}-1} \log\left( - \nabla^2 \right) \\
    + \alpha_0 c_0^y \left( 2 \log(c_0) \tan \left( \frac{\pi y}{2} \right) + \pi \mathrm{sec}^2 \left( \frac{\pi y}{2} \right) \right) \left(- \nabla^2 \right)^{\tfrac{y+1}{2}-1} + \frac{\eta}{2} \left(- \nabla^2 \right)^{\tfrac{y+1}{2}-1} \log\left( - \nabla^2 \right)
\end{multline}

\section{$k$-Space Pseudospectral Time Domain Implementation of Time-Reversal-Based Gradient Computation} \label{sec:kSpaceImp}

Here we sketch how to implement the proposed memory efficient gradient computation \eqref{eq:TRBwdWave} using a $k$-space TD scheme to solve the underlying first order system of equations \eqref{eq:MonCon}-\eqref{eq:PreDenRel}. It follows the same approach implemented in \kwave/ \cite{TrCo10}. For a more detailed description, we refer to the \kwave/ manual\footnote{{\textcolor{darkblue}{\texttt{www.k-wave.org}}}}. The $k$-space TD method is a collocation scheme that interpolates between the collocation points using a truncated Fourier series. This allows gradients to be calculated using the \termabb{fast Fourier transform}{FFT} which leads to fast implementations on GPU architectures. The `$k$-space' in the name refers to a correction $\kappa$ applied in the Fourier domain to account for the finite difference approximation of the time derivative. This correction is exact for a homogeneous medium \cite{TaMaWa02} and reduces the errors for acoustically heterogeneous media. The spatial domain $\Omega$ has to be embedded into a rectangular region, which is then discretised by a regular grid of $N$ collocation points with grid size $\Delta x$. To mimic free-space propagation, a \termabb{Perfectly Matched Layer}{PML} absorbing boundary is wrapped around the box to damp outgoing waves without reflecting them. In this work, the size of the PML layer is automatically chosen to be between $10-20$ voxels by minimizing the the highest prime factor of the total grid size (cf. function \texttt{getOptimalPMLSize.m} in \kwave/). This will lead to efficient fft executions.

We discretize the measurement time interval $[0,T]$ by $t_n =  n \Delta t$, $n=0,\ldots,N_t$, $\Delta t = T/N_t$. The discrete pressure is denoted by $p$ and the $\xi$-component of the discrete particle velocity by $v_\xi$ ($\xi$ stands for one of the spatial components $x, y, z$ of these vector fields). As the terms containing $\nabla \rho_0(x)$ in \eqref{eq:MonCon}-\eqref{eq:PreDenRel} cancel out, they are omitted \cite{TrJaReCo12}. The scalar density $\rho$ is split into three parts, $\rho_\xi$. This is nonphysical but allows anisotropic absorption to be included within the PML. The \emph{k-space derivative} and the \emph{k-space operator} $\kappa$ are defined as 
\begin{equation}
\frac{\partial}{\partial \xi} g \mydef \cF^{-1} \left\lbrace i k_\xi \kappa \cF \left\lbrace g \right\rbrace \right\rbrace; \qquad \kappa = \text{sinc}\left(\frac{c_{{\mathrm ref}} \Delta t}{2} \sqrt{k_x^2+k_y^2+k_z^2} \right), \label{eq:DisPartDev} 
\end{equation}
where $k_\xi \in \R^N$ is the discrete wavevector in $\xi$ direction and all multiplications between $N$-dimensional vectors are understood as componentwise. $c_{\mathrm ref}$ is a reference sound speed, chosen to ensure stability. For consistency during the FWI inversion, we fix $c_{\mathrm ref} = c^{max}_0$. Grid-staggering is incorporated into the calculation of the gradients. This means that a spatial translation of $\pm \Delta \xi / 2$ is introduced as
\begin{equation}
\frac{\partial^\pm}{\partial \xi} g \mydef \cF^{-1} \left\lbrace i k_\xi \kappa e^{\pm i k_\xi \Delta \xi/2} \cF \left\lbrace g \right\rbrace \right\rbrace. \label{eq:DisPartDevSt} 
\end{equation}
Staggering in time will be included by interleaving the gradient and updates steps. Multiplication operators $\Lambda_\xi$ and $\Lambda^s_\xi$ implement the PML on the normal and staggered grid, respectively. The measurement operator $M$, i.e., the receiving US transducers, is modeled by a linear mapping between spatial grid and measurement channels and is also denoted by $M$. This assumes that it acts instantaneously, i.e., in iteration $n$, we approximate $f(t_n) = M p$ (otherwise, the pressure fields at the sensor locations need to be buffered). In the same way, a source operator $S$ will embed the discrete pressure source $s(t_n)$ into the spatial grid. After this, $S s(t_n)$ is divided by $c_0$ pointwise and scaled by $2 / (3 \Delta x)$ to convert it into a mass source, cf. Section 2.4 in the \kwave/ manual.  To implement the TR-based gradient \eqref{eq:TRBwdWave}, we need to extract the pressure on the boundary of sample, $\partial \Gamma$, during the forward wave simulation. In our collocation scheme, we will simply extract the pressure in a layer of grid points with a defined thickness as examined in Section \ref{subsec:ValTRGrad}, which we will denote as $p(\partial \Gamma)$ for simplicity. The scheme to compute the TR-based gradient consists of three wave simulations, \textit{forward}, \textit{adjoint}, and \textit{time-reversal} of which the latter two are interleaved. For simplicity, we only consider $u = c_0$ here (as in the numerical studies) and any scaling factors like $\Delta t$ have been omitted if they cancel out in the end. Algorithm \ref{algo:TRImplementation} describes the whole iteration. Step \ref{step:cumsum} corrects for the fact that the first order scheme effectively uses the time derivative of the mass source term provided, not the term itself, cf.  \eqref{eq:WaveSecondOrder}. Code will be released as part of a more general \textsc{Matlab} toolbox for ultrasonic and photoacoustic image reconstruction.

\begin{algorithm}
    \caption{Time Reversal Based Speed-Of-Sound Gradient Computation}\label{algo:TRImplementation}
    \begin{algorithmic}[1]
        \State \textbf{input:} $s$
        \Statex
        \Statex \textbf{forward iteration}
        \State $p, v_\xi, \rho_\xi \gets 0$ for all $\xi \in [x,y,z]$ \Comment{initilization forward iteration} 
        \For{$n \, = \, 0, \ldots, N_t$} \Comment{$t_0 = 0$ and $t_{N_t} = T$}
            \State $v_\xi \hspace{0.1em} \gets \Lambda^s_\xi \left( \Lambda^s_\xi v_\xi - \frac{\Delta t}{\tilde{\rho}_0} \frac{\partial^+}{\partial \xi} p \right)$ \Comment{momentum conservation, cf. \eqref{eq:MonCon}}
            \State $\rho_\xi\hspace{0.1em} \gets \Lambda_\xi \left( \Lambda_\xi \rho_\xi - \Delta t \tilde{\rho}_0 \frac{\partial^-}{\partial \xi} v_\xi \right) + \frac{2 \Delta t}{3 \Delta x c_0} \,  S s(t_n)$ \Comment{mass conservation, cf. \eqref{eq:MassCon}}
            \State $p\hspace{0.6em} \gets c_0^2 \left(\rho_x + \rho_y + \rho_z \right)$  \Comment{pressure-density relation, cf. \eqref{eq:PreDenRel}}
            \State $f^{n} \gets M p$ \Comment{simulate measurement}
            \State $g^{n}\hspace{0.1em} \gets p(\partial \Gamma)$  \Comment{extract boundary data}
        \EndFor
        \Statex
        \Statex \textbf{setup interleaved adjoint and time reversal iterations}
        \State $\nabla \cD, p, v_\xi, \rho_\xi, p^\triangleleft, v^\triangleleft_\xi, \rho^\triangleleft_\xi \gets 0$ for all $\xi \in [x,y,z]$ \Comment{initilization adjoint and TR iterations} 
        \State $q^n \gets f^{N_t-n} - f^\delta(T - t_n)$ for all $n$ \Comment{time-reversed residual as adjoint source}
        \State $\cD \gets \tfrac{1}{2} \sqnorm{q}$ \Comment{compute loss}
        \State $q \gets \texttt{cumsum}(q)$ \label{step:cumsum} \Comment{integrate $q$ in time}
        \State $p^\triangleleft(\partial \Gamma) \gets g^{N_t}$ \Comment{initial value}
        \State $g^n \gets g^{n} / (3 \, c_0(\partial \Gamma)^2)$ for all $n$ \Comment{scale from pressure to density fields} 
        \State $p_{rec}^0 \gets p^\triangleleft(\Gamma)$ \Comment{reconstructed forward field in $\Gamma$ at the current time step}
        \Statex
        \Statex \textbf{advance time reversal iteration by one step to get ahead of adjoint iteration:}
        \State $v^\triangleleft_\xi \gets  \Lambda^s_\xi \left( \Lambda^s_\xi v^\triangleleft_\xi - \frac{\Delta t}{\tilde{\rho}_0} \frac{\partial^+}{\partial \xi} p^\triangleleft \right)$ 
        \State $\rho^\triangleleft_\xi \gets  \Lambda_\xi \left( \Lambda_\xi \rho^\triangleleft_\xi - \Delta t \tilde{\rho}_0 \frac{\partial^-}{\partial \xi} v^\triangleleft_\xi \right)$ 
        \State $\rho^\triangleleft_\xi(\partial \Gamma) \gets g^{N_t-1}$ \Comment{enforce Dirichlet condition} 
        \State $p^\triangleleft \gets c_0^2 \left(\rho^\triangleleft_x + \rho^\triangleleft_y + \rho^\triangleleft_z \right)$ 
        \State $p_{rec}^+ \gets p^\triangleleft(\Gamma)$ \Comment{reconstructed forward field at $n+1$}
        \Statex
        \Statex \textbf{advance adjoint and time reversal iterations simultaneously:}
        \For{$n \, = \, 1, \ldots, N_t-1$} 
            \State $v_\xi \gets  \Lambda^s_\xi \left( \Lambda^s_\xi v_\xi - \frac{\Delta t}{\tilde{\rho}_0} \frac{\partial^+}{\partial \xi} p \right)$
            \State $\rho_\xi \gets  \Lambda_\xi \left( \Lambda_\xi \rho_\xi - \Delta t \tilde{\rho}_0 \frac{\partial^-}{\partial \xi} v_\xi \right) + \frac{2 \Delta t}{3 \Delta x c_0} \,  M^T q^n$
            \State $p \gets \quad c_0^2 \left(\rho_x + \rho_y + \rho_z \right)$
            \State $v^\triangleleft_\xi \gets  \Lambda^s_\xi \left( \Lambda^s_\xi v^\triangleleft_\xi - \frac{\Delta t}{\tilde{\rho}_0} \frac{\partial^+}{\partial \xi} p^\triangleleft \right)$ 
            \State $\rho^\triangleleft_\xi \gets  \Lambda_\xi \left( \Lambda_\xi \rho^\triangleleft_\xi - \Delta t \tilde{\rho}_0 \frac{\partial^-}{\partial \xi} v^\triangleleft_\xi \right)$ 
            \State $\rho^\triangleleft_\xi(\partial \Gamma) \gets g^{N_t-n-1}$ 
            \State $p^\triangleleft \gets c_0^2 \left(\rho^\triangleleft_x + \rho^\triangleleft_y + \rho^\triangleleft_z \right)$ 
            \State $p_{rec}^- \gets p_{rec}^0$
            \State $p_{rec}^0 \gets p_{rec}^+$
            \State $p_{rec}^+ \gets p^\triangleleft(\Gamma)$ 
            \State $\nabla \cD \gets \nabla \cD +  \frac{2}{\Delta t c_0^3} \left(p_{rec}^+ - 2 p_{rec}^0 + p_{rec}^-\right) p(\Gamma)$ \Comment{assemble gradient, cf. \eqref{eq:Gradc0}, $p$ corresponds to $q^*$}
        \EndFor
        \Statex 
        \State \textbf{return:} $\cD, \nabla \cD$
    \end{algorithmic}
\end{algorithm}

\section{Multigrid for $k$-Space Pseudospectral Schemes} \label{sec:MultiGridDetails}

\begin{figure}[tb!]
   \centering
\includegraphics[width=\textwidth]{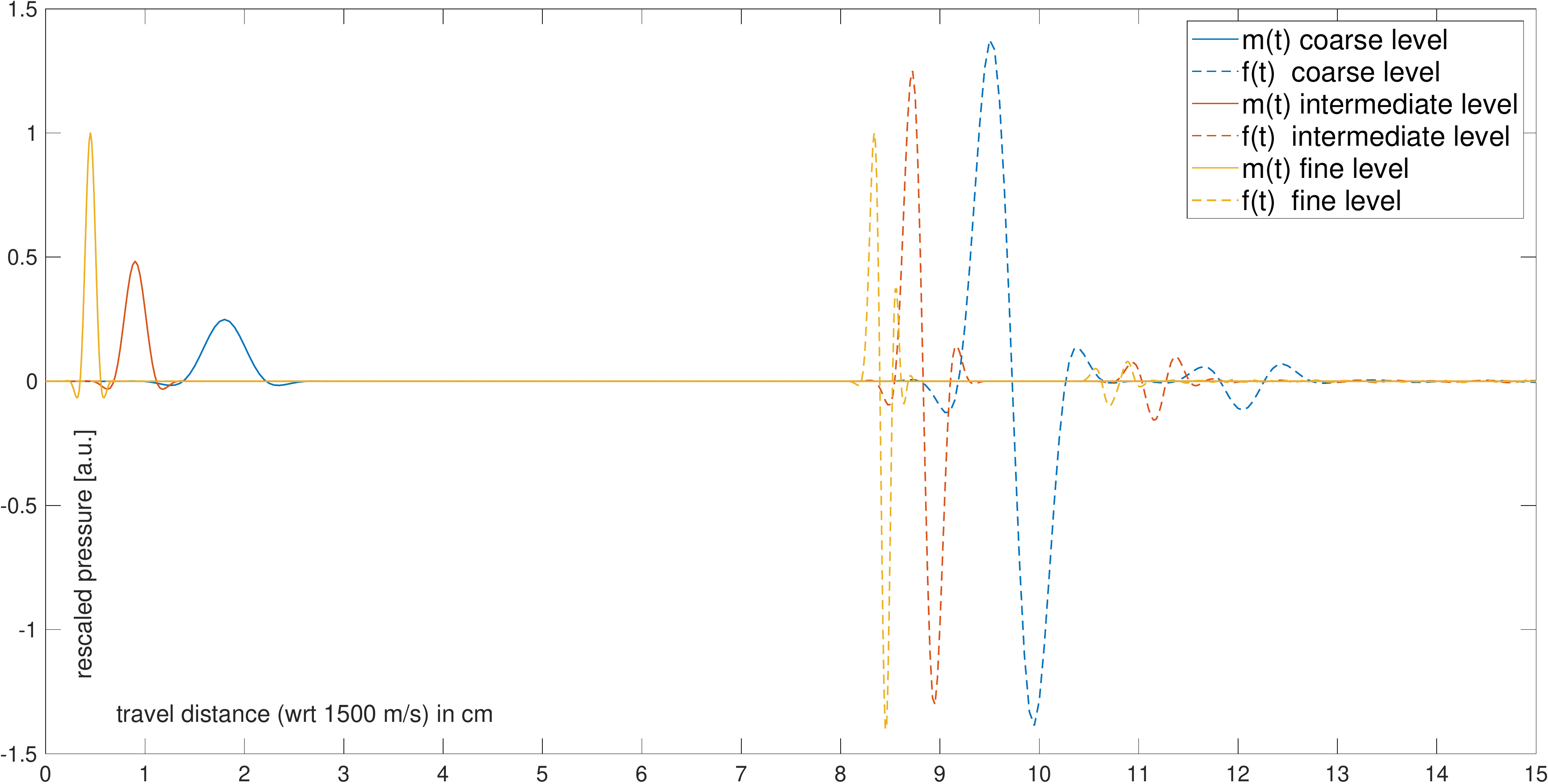}
\caption{Mass source signals $m(t)$ and corresponding simulated measurements $f(t)$ on different spatio-temporal grids ($\eta = 2$) for a random source-receiver pair. The signals on coarser grids result from Fourier interpolation of the fine grid signals followed by low-pass filtering. The signals have been rescaled for comparison. On each level, the temporal filtering used to remove the highest frequencies not supported by the grid introduces a different phase shift. However, as the shift is the same for source and measured time series, it does not affect the FWI results.}
   \label{fig:MultiGridSignals}
\end{figure}

In the Fourier collocation scheme described above, the discrete, distributed field parameters $(p, v, \rho)$ implicitly represent band-limited interpolants of the corresponding infinite dimensional variables (cf. Section 2.8. in the \kwave/ manual). As such, Fourier/trigonometric interpolation is the natural way to transfer variables from one spatial grid into another. For regular spatial grids, this can be implemented efficiently using FFTs (cf. function \texttt{interpftn.m} in \kwave/), and can be combined with smoothing in Fourier space. To implement the coarse-to-fine initialization strategy described in Section \ref{subsec:MultGrid} and examined in Section \ref{subsec:ValMultGrid}, we first need to down-sample the temporal dimension of data $f^{\delta}$ and our model of the corresponding (mass) source $m$ to a coarser temporal grid. This is done via a combination of Fourier interpolation and low-pass filtering with a Kaiser window (cf. function \texttt{filterTimeSeries.m} in \kwave/) to remove frequencies not supported by the numerical scheme. The data additionally needs to be scaled by $\eta^{d-1}$. Figure~\ref{fig:MultiGridSignals} illustrates the results. After the FWI solution $u^c$ is computed on the coarse grid, we need to up-sample it to initialize the FWI on the next, finer grid. For this, we use a combination of Fourier interpolation and smoothing with a Blackman window  (cf. Section 2.8. in the \kwave/ manual): $u^f \gets \texttt{interpftn(}u^c, [N_x^f, N_y^f, N_z^f], \texttt{\textquotesingle Blackman\textquotesingle)}$, where $[N_x^f, N_y^f, N_z^f]$ are the dimensions of the fine spatial grid. 

\section{L-BFGS with Stochastic Gradient Estimates} \label{sec:SLBFGS}

We want to minimize an energy $\cJ(u)$ as in \eqref{eq:FWI} while only having access to an unbiased estimator $g(u, \theta, n)$ of $(\cJ(u), \nabla \cJ(u))$. Here, $\theta \in [0,1]$ determines the approximation strength of $g$ where $\theta = 0$ corresponds to no approximation, e.g., $g(u,0,\cdot) = \cJ(u)$, and $\theta = 1$ to maximal approximation, e.g., by using a data sub-set $\cS$ in $g_{\cS}$ of size one. In our implementation, $\theta$ can be a function of the iteration count $i$. The last input $n \in \mathbb{N}$ in $g(u, \theta, n)$ controls the the stochastic approximation performed by $g$. For instance, it could seed the random generator used to select $\cS$ in $g_\cS$. Important for the algorithm is only that it leads to the same type of approximation for different $u$, e.g., for a given $\theta, n$, $g_{\cS}(\cdot, \theta, n)$ always uses the same subset of data $\cS$. 

Algorithm \ref{algo:SLBFGS} describes the stochastic L-BFGS (SLBFGS) scheme used in this work, which is based on the one presented in \cite{FaGlGi17}. See Section 6.2. in \cite{BoCuNo18} and references therein for a detailed discussion of stochastic quasi-Newton schemes. The key idea is that for computing the vectors used to construct the low-rank approximation of the inverse Hessian, $H_k$, only iterate pairs evaluated with the same stochastic approximation should be used. This increases the number of gradient computations to at least two per iteration (more if line-search is used) but the second gradient evaluation can be utilized to compute a second update. In \cite{FaGlGi17}, no bound constraints $u_{min} \leq u \leq u_{max}$ were considered. Here, we include them but assume that they are not active at the vicinity of the solution. This means that they do not stabilize the image reconstruction problem. Rather, they merely safeguard against running the numerical wave simulation scheme with values of the parameter fields $u$ that would lead to instability.
For this reason, we enforce them by a simple projection of the updated variable $u^k$ at the end of each iteration (step \ref{setp:Proj}). We use progressive iterate averaging with $w(i) = i^3$, but only start averaging after the estimated energy increases for the first time (step \ref{setp:AvgAct}). We use history size $m_h = 64$ throughout the paper. Code will be released as part of a more general \textsc{Matlab} toolbox for deterministic and stochastic reconstruction.

\begin{algorithm}
    \caption{Stochastic Limited Memory BFGS (SLBFGS)}\label{algo:SLBFGS}
    \begin{algorithmic}[1]
        \State \textbf{inputs}: $u, H_0, w, m_h, N_{evl}, \nu, \texttt{do\_linesearch}$ \Comment{initial value, initial approx of inverse Hessian, weight function, history size, maximal function evaluations, step size, boolean determining if linesearch is used} 
        \Statex
        \State $\bar{u}, \bar{w} \gets 0$ \Comment{averaged iterates, accumulated weights}
        \State $J_0 \gets \infty$ \Comment{energy estimate}
        \State $S, Y \gets [\;]$            \Comment{initialize buffer to empty array}
        \State \texttt{average\_iterates} $\gets$ \texttt{false} \Comment{iterate averaging is inactive}
        \State $i \gets 0$            \Comment{number of iterations}
        \Statex
        \While{$n_{evl} < N_{evl} $} \Comment{iterate until max number of function evaluations reached}
            \State $i \gets i + 1$   \Comment{increase iteration count}
            \State $F_u, G_u \gets g(u, \theta(i), i)$ \Comment{evaluate estimator}
            \State $n_{evl} \gets n_{evl} + 1$   \Comment{update function evaluation count}
            \State $z \gets - \texttt{twoLoopRecursion}(G_u, S, Y, H_0)$ \Comment{implements $z = - H_k G_u$}
            
            \If{\texttt{do\_linesearch}}
                \State $F_z, G_z, \nu, m_{evl} \gets \texttt{lineSearch}(u, z, g(\,\cdot\,, \theta(i), i), \nu)$ \Comment{line-search on $f(u) := g(u,\theta(i),i)$}     
            \Else
                \State $F_z, G_z \gets g(u + \nu z, \theta(i), i)$ \Comment{evaluate estimator at $u + \nu z$}
                \State $m_{evl} \gets 1$
            \EndIf
            \State $n_{evl} \gets n_{evl} + m_{evl}$   \Comment{update function evaluation count}
            \State $S \gets \texttt{updateBuffer}(\nu z, S, m_h)$     \Comment{update buffer $S$}
            \State $Y \gets \texttt{updateBuffer}(G_z - G_u, Y, m_h)$ \Comment{update buffer $Y$}
            \State $z \gets z - \texttt{twoLoopRecursion}(G_z, S, Y, H_0)$ \Comment{implements $z = z - H_k G_z$}
            \State $u \gets u + \nu z$  \Comment{update $u$}
            \State $u \gets \texttt{max}( \texttt{min}(u, u_{max}), u_{min})$  \label{setp:Proj} \Comment{project onto bound constraints}
            \State $J_i \gets \texttt{min}(F_u, F_z)$ \Comment{update energy estimate}
            \If{\texttt{average\_iterates}}
                \State $\bar{u} \gets \bar{u} + w(i) u$  \Comment{update accumulated weighted iterates}
                \State $\bar{w} \gets \bar{w} + w(i)$  \Comment{update accumulated weights}
                \State $\hat{u} \gets \bar{u} / \bar{w}$ \Comment{update estimate for minimizer}
            \Else
                \If{$J_i > J_{i-1}$} \Comment{first increase in estimated energy activated averaging}
                    \State $\texttt{average\_iterates} \gets \texttt{true}$ \label{setp:AvgAct}             
                \EndIf
                \State $\hat{u} \gets u$ \Comment{update estimate for minimizer}
            \EndIf
        \EndWhile
        \Statex 
        \State \textbf{return:} $\hat{u}, J$
    \end{algorithmic}
\end{algorithm}

\section{Delayed Source Encoding for Time-Invariant Systems} \label{sec:SourceEncodingTimeInvariantSys}

\begin{figure}[tb!]
   \centering
\includegraphics[width=\textwidth]{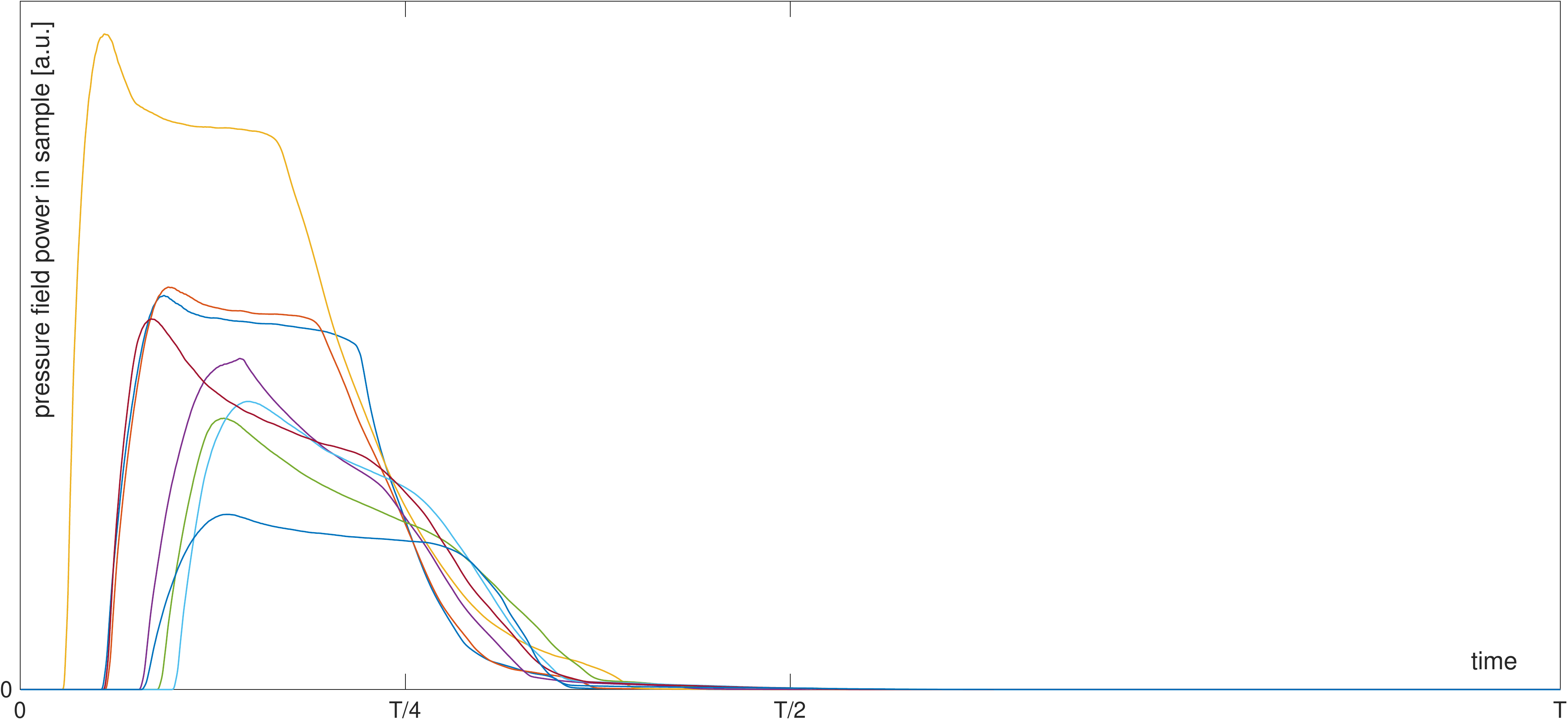}
\caption{Temporal evolution of the pressure field power $P_i(t)$ for 8 different sources.}
   \label{fig:SignalDecay}
\end{figure}

Here, we examine the conventional and delayed source encoding described in Section \ref{subsec:SGDandSE} in more detail. It is crucial that we assumed that the acoustic properties are known everywhere outside the sample $\Gamma$ (e.g., $\mathrm{supp}(u) \subset \Gamma$), and that $p(x,t)$ is not practically measurable for $t > T$ anymore. Let $\hat p$ and $\hat q^*$ be the forward and adjoint wave fields corresponding to the delayed encoded sources $\hat s$ and measurements $\hat f^\delta$, cf. \eqref{eq:DelaySE}. By the linearity and time invariance of the wave equation, we have $\hat p(x,t) = \sum_i w_i p_i(x, t - d_i \tau)$ and $\hat q^*(x,t) = \sum_i w_i q_i^*(x, t- d_i \tau)$. With this, we get (cf. Section \ref{subsec:FWI} and \eqref{eq:Gradc0}): 
\begin{align}
\nabla_{u} \cD \left( \hat{f}(u) , \hat{f}^{\delta} \right) =    \int_0^{T+\tau}  \left( \frac{\partial A}{\partial u} \hat p(x,t) \right) \hat q^*(x, t) \, dt=    \sum_{i,j}^{n_s}  w_i w_j \int_0^{T+\tau} \left(\frac{\partial A}{\partial u}  p_i(x,t - d_i \tau) \right)  q_j^*(x, t - d_j \tau) \, dt \nonumber \\
=   \sum_i^{n_s} w_i^2   \int_0^{T + \tau}  \left(\frac{\partial A}{\partial u}  p_i(x,t - d_i \tau) \right)  q_i^*(x, t - d_i \tau)  \, dt +  \sum_{i \neq j}^{n_s} w_i w_j \int_0^{T+\tau}  \left(\frac{\partial A}{\partial u}  p_i(x,t - d_i \tau) \right)  q_j^*(x, t - d_j \tau) \, dt  \nonumber \\
=   \sum_i^{n_s} w_i^2   \int_0^{T}  \left(\frac{\partial A}{\partial u}  p_i(x,t) \right)  q_i^*(x, t) \, dt  +  \sum_{i \neq j}^{n_s} w_i w_j \int_0^{T+\tau}  \left(\frac{\partial A}{\partial u}  p_i(x,t - d_i \tau) \right)  q_j^*(x, t - d_j \tau) \, dt  \nonumber \\
=   \sum_i^{n_s} w_i^2 \nabla_{u} \cD \left( f_i(u) , f_i^\delta \right)    +  \sum_{i \neq j}^{n_s} w_i w_j \int_0^{T+\tau}  \left(\frac{\partial A}{\partial u}  p_i(x,t - d_i \tau) \right)  q_j^*(x, t - d_j \tau) \, dt  
\end{align}
With $\Exp \left[ w \right] = 0$, $\Cov[w] = I$, it is easy to see now that $\Exp \left[ \nabla_{u} \cD \left( \hat{f}(u) , \hat{f}^{\delta} \right) \right] = \nabla \cJ(u)$. If we choose the Rademacher distribution ($w_i = \pm 1$ with equal probability), we have that $w_i^2 = 1$ and
\begin{equation}
\nabla_{u} \cD \left( \hat{f}(u) , \hat{f}^{\delta} \right) = \nabla \cJ(u) +  \sum_{i \neq j}^{n_s} w_{i} w_{j} X_{ij}, \quad \text{where} \quad X_{ij} := \int_0^{T+\tau}  \left(\frac{\partial A}{\partial u}  p_i(x,t - d_i \tau) \right)  q_j^*(x, t - d_j \tau)  \, dt     
\end{equation}
is the cross talk between sources $i$ and $j$. To examine the impact of delays ($\tau > 0$) vs. conventional source encoding ($\tau = 0$), we assume $\bar{d} = d_j - d_i > 0$ without loss of generality and shift the time in the integral by $d_i \tau$ to get 
\begin{equation}
X_{ij}(\tau) = \int_0^{T+\tau}  \left(\frac{\partial A}{\partial u}  p_i(x,t) \right)  q_j^*(x, t - \bar{d} \tau)  \, dt 
\end{equation}
For $\tau \geq T / \bar{d}$, the delay between both sources is bigger than $T$, and from our assumptions, it follows that $X_{i,j}(\tau) = 0$. To obtain generic estimates for $\tau < T / \bar{d}$, one needs to estimate how fast the field power 
\begin{equation}
    P_i(t) := \int_\Gamma p_i^2(x,t) \, dx
\end{equation}
decays inside $\Gamma$. For such decay estimates, see Section ``Local Energy Decay Estimates'' in \cite{KuKu11} and references therein. In 3D, there is a $T_d$ such that $P_i(t) < a e^{- b t}$ with $a, b > 0$ for all $t > T_d$. In Figure~\ref{fig:SignalDecay}, we plot $P_i(t)$ for 8 different sources in the scenario used in the numerical studies to show that practically, the decay is quite rapid, which is also confirmed by the results in Section \ref{subsec:ValStochOpt}, cf. Figure~\ref{fig:StochGradEst}. 

\section{Noise Sensitivity} \label{sec:Noise}

\begin{figure}[tb!]
   \centering
\hfill
\subfloat[][\label{subfig:NoiseTimeTrace}]{\includegraphics[height=0.4\textwidth]{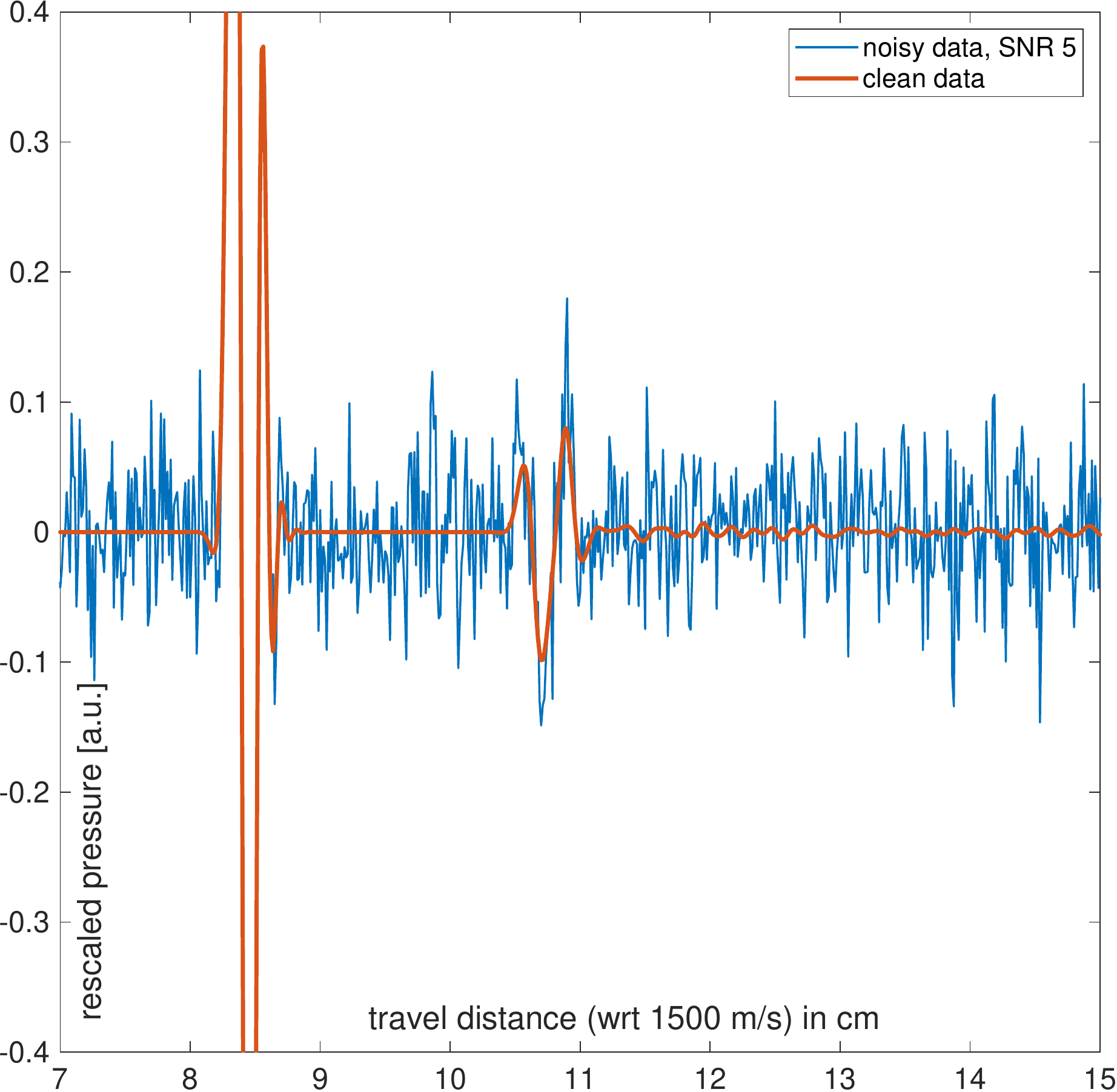}}
\hfill
\subfloat[][\label{subfig:NoiseLogErr}]{\includegraphics[height=0.4\textwidth]{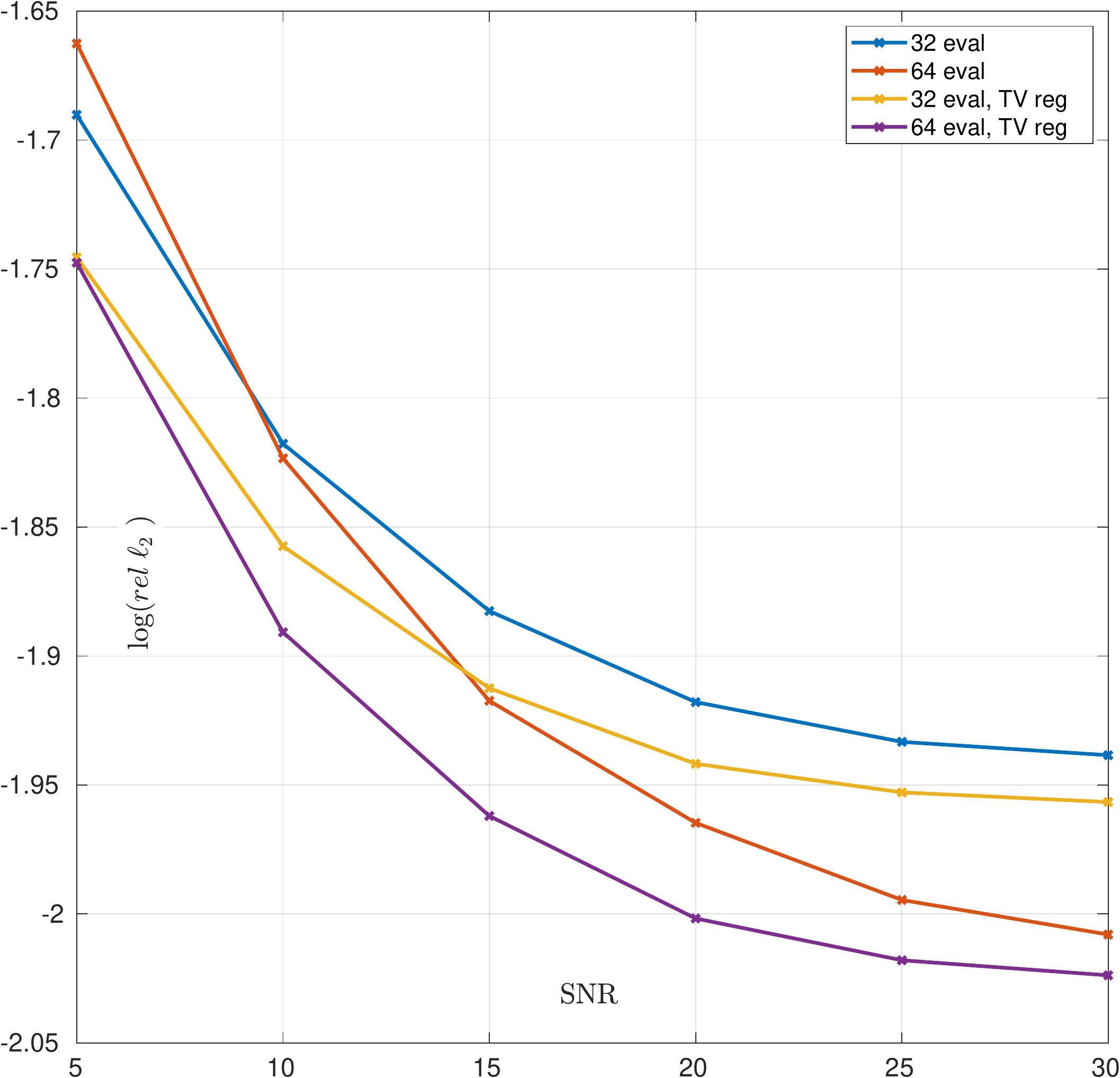}}
\hfill
\caption{\protect\subref{subfig:NoiseTimeTrace} Comparison of clean data $f(t)$ and noisy data $f(t) + \varepsilon(t)$ for SNR 5: The clean data is the same time trace shown as "$f(t)$ fine level" in Figure \ref{fig:MultiGridSignals}. The noise signal $\varepsilon(t)$ for SNR 10,15,20,25,30 is the same just multiplied with $10^{\frac{5- SNR}{20}}$ = 0.56, 0.32, 0.18 , 0.10, 0.056. \protect\subref{subfig:NoiseLogErr} Quantitative results of the noise sensitivity study.}
   \label{fig:NoiseLogErr}
\end{figure}

In this section, we demonstrate the impact of adding measurement noise to the simulated data. The noise characteristics of UST systems can differ quite a lot from conventional linear array US systems, e.g., both the US transducers and the data acquisition (DAQ) system of the PAMMOTH scanner were custom designed and a detailed noise characterization will be included in forthcoming work. Here, we will assume that an additive Gaussian noise model is adequate to summarize the different noise contributions and further that a decorrelation transform will render its covariance matrix sufficiently close to a scaled identity: 
\begin{equation}
f^{\delta} = A(u^\dagger) + \varepsilon, \qquad \text{with} \qquad \varepsilon \sim \cN(0, \sigma^2 I),
\end{equation}
In the following, we will use values of $\sigma$ determined by the \textit{signal-to-noise (SNR)} level defined in terms of the \textit{root mean square (RMS)} of the clean signal $f$ as 
\begin{equation}
\sigma = RMS(f) \, 10^{-SNR/20}.
\end{equation} 
Figure \ref{subfig:NoiseTimeTrace} shows an example of a clean and noisy time trace for SNR $5$. \\
Without measurement noise ($\varepsilon = 0$), all noise-like features in the reconstructed images were a result of using stochastic gradient estimators in the optimization schemes to solve the FWI problem \eqref{eq:FWI}, i.e., a purely numerical phenomena. This "gradient-noise" vanishes as the iteration progresses (cf. Figure \ref{fig:PreStudy1}) or if more accurate gradient estimators are used (cf. Figure \ref{fig:FullDemo}). From the way that $\nabla \cJ(u)$ is computed as presented in Section \ref{subsec:FWI} is is clear that any measurement noise $\varepsilon > 0$ will be back-propagated into the image alongside with the data residual. As the iteration progresses, the norm of this type of image noise will grow while its spectral content changes. This phenomena arises from the fact that in principle, we use a Landweber-type iteration to solve an ill-conditioned inverse problem \eqref{eq:USTmodel}. See \cite{Vo02,BuKaNe11} for a general introduction into this topic. 
Here, we illustrate these two opposing effects by first using exactly the same settings as in Section \ref{sec:FullDemo} on a single GPU for different SNR values. The first column of Figures \ref{fig:NoiseStudy}, \ref{fig:NoiseStudyDiff} show that while SNR 30 is visually indistinguishable from clean data result (Figure \ref{subfig:RecDesktop}) the image noise increases towards SNR 5. Then we change the number of gradient evaluations on finest spatial level from 32 to 64. The second column of Figures \ref{fig:NoiseStudy}, \ref{fig:NoiseStudyDiff} and the first two plots in Figure \ref{subfig:NoiseLogErr} show the results and confirm that for the highest noise level (SNR = 5), doing more iterations now indeed leads to a worse result. To prevent this, we will need to add a regularization functional to the data fidelty $\sum_i^{n_s} \cD_i(u)$  in \eqref{eq:FWI}. For simplicity, we use a smoothed version of the popular \textit{Total Variation (TV)} functional as described in \cite{Vo02} here: 
\begin{equation}
\min_{u \in \cU} \, \sum_i^{n_s} \cD_i(u) + \alpha \, \text{TV}(u)\qquad \text{with} \qquad \text{TV}(u) := \sum_i \sqrt{\left(\partial_x u \right)^2_i + \left(\partial_y u \right)^2_i + \left(\partial_z u \right)^2_i + \beta^2 } - \beta,   \label{eq:RegFWI}
\end{equation}
where we set $\beta$ to \SI{1}{\meter / \second} and fix $\alpha > 0$ to an illustrative value that ensures observable impact of the regularization (for optimal results, one would need to choose it in dependence of $\sigma$). The third and forth columns of Figures \ref{fig:NoiseStudy}, \ref{fig:NoiseStudyDiff} and the second two plots in Figure \ref{subfig:NoiseLogErr} show that now the errors do not increase with increasing iteration count. Note however that it is not straight forward to predict the overall impact of adding regularization to FWI and in particular the question which type of functionals are advantageous for breast UST needs further research.

\newcommand{\heightNoise}{0.3}
\begin{figure}[tb!]
\centering
\subfloat[][SNR 30, 32 eval\label{subfig:SNR30I32}]{\includegraphics[height=\heightNoise\textwidth]{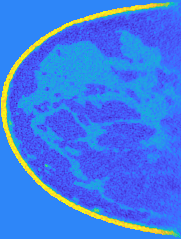}}
\hfill
\subfloat[][SNR 30, 64 eval\label{subfig:SNR30I64}]{\includegraphics[height=\heightNoise\textwidth]{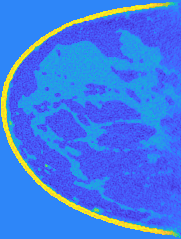}}
\hfill
\subfloat[][SNR 30, TV reg, 32 eval\label{subfig:SNR30I32TV}]{\includegraphics[height=\heightNoise\textwidth]{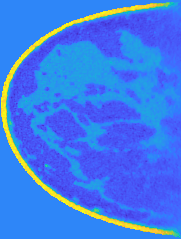}}
\hfill
\subfloat[][SNR 30, TV reg, 64 eval\label{subfig:SNR30I64TV}]{\includegraphics[height=\heightNoise\textwidth]{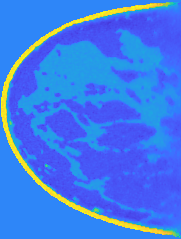}}\\
\subfloat[][SNR 20, 32 eval\label{subfig:SNR20I32}]{\includegraphics[height=\heightNoise\textwidth]{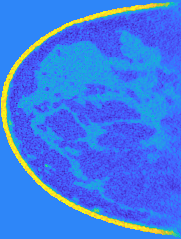}}
\hfill
\subfloat[][SNR 20, 64 eval\label{subfig:SNR20I64}]{\includegraphics[height=\heightNoise\textwidth]{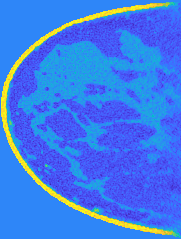}}
\hfill
\subfloat[][SNR 20, TV reg, 32 eval\label{subfig:SNR20I32TV}]{\includegraphics[height=\heightNoise\textwidth]{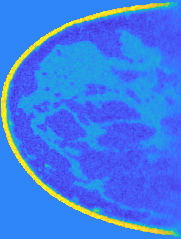}}
\hfill
\subfloat[][SNR 20, TV reg, 64 eval\label{subfig:SNR20I64TV}]{\includegraphics[height=\heightNoise\textwidth]{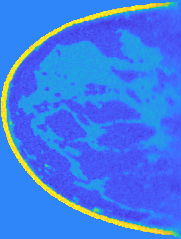}}\\
\subfloat[][SNR 10, 32 eval\label{subfig:SNR10I32}]{\includegraphics[height=\heightNoise\textwidth]{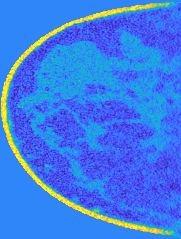}}
\hfill
\subfloat[][SNR 10, 64 eval\label{subfig:SNR10I64}]{\includegraphics[height=\heightNoise\textwidth]{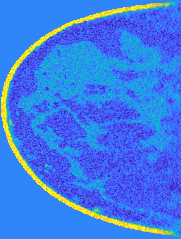}}
\hfill
\subfloat[][SNR 10, TV reg, 32 eval\label{subfig:SNR10I32TV}]{\includegraphics[height=\heightNoise\textwidth]{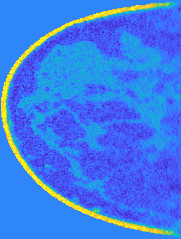}}
\hfill
\subfloat[][SNR 10, TV reg, 64 eval\label{subfig:SNR10I64TV}]{\includegraphics[height=\heightNoise\textwidth]{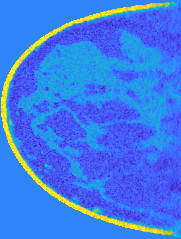}}\\
\subfloat[][SNR 5, 32 eval\label{subfig:SNR5I32}]{\includegraphics[height=\heightNoise\textwidth]{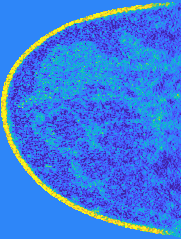}}
\hfill
\subfloat[][SNR 5, 64 eval\label{subfig:SNR5I64}]{\includegraphics[height=\heightNoise\textwidth]{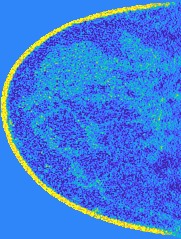}}
\hfill
\subfloat[][SNR 5, TV reg, 32 eval\label{subfig:SNR5I32TV}]{\includegraphics[height=\heightNoise\textwidth]{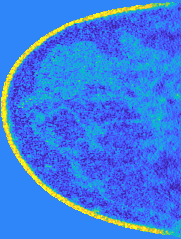}}
\hfill
\subfloat[][SNR 5, TV reg, 64 eval\label{subfig:SNR5I64TV}]{\includegraphics[height=\heightNoise\textwidth]{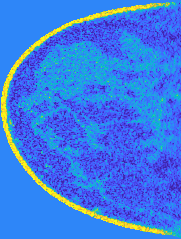}}\\
\caption{Qualitative results of the noise sensitivity study. Shown are different combinations of noise level (SNR), maximal number of gradient evaluations (32 eval vs 64 eval) and of using additional Total Variation regularization.}
   \label{fig:NoiseStudy}
\end{figure}

\begin{figure}[tb!]
\centering
\subfloat[][SNR 30, 32 eval\label{subfig:SNR30I32Diff}]{\fbox{\includegraphics[height=\heightNoise\textwidth]{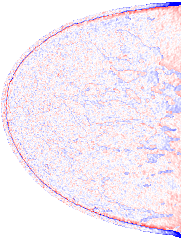}}}
\hfill
\subfloat[][SNR 30, 64 eval\label{subfig:SNR30I64Diff}]{\fbox{\includegraphics[height=\heightNoise\textwidth]{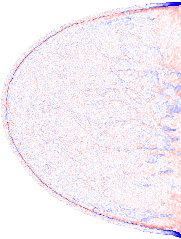}}}
\hfill
\subfloat[][SNR 30, TV reg, 32 eval\label{subfig:SNR30I32TVDiff}]{\fbox{\includegraphics[height=\heightNoise\textwidth]{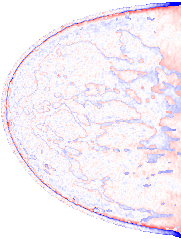}}}
\hfill
\subfloat[][SNR 30, TV reg, 64 eval\label{subfig:SNR30I64TVDiff}]{\fbox{\includegraphics[height=\heightNoise\textwidth]{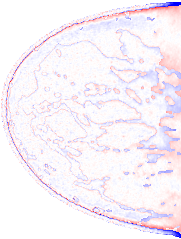}}}\\
\subfloat[][SNR 20, 32 eval\label{subfig:SNR20I32Diff}]{\fbox{\includegraphics[height=\heightNoise\textwidth]{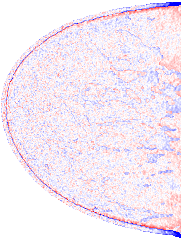}}}
\hfill
\subfloat[][SNR 20, 64 eval\label{subfig:SNR20I64Diff}]{\fbox{\includegraphics[height=\heightNoise\textwidth]{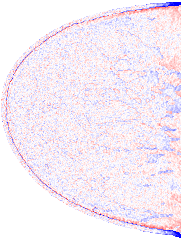}}}
\hfill
\subfloat[][SNR 20, TV reg, 32 eval\label{subfig:SNR20I32TVDiff}]{\fbox{\includegraphics[height=\heightNoise\textwidth]{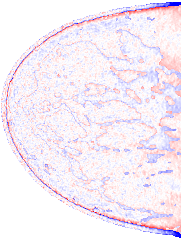}}}
\hfill
\subfloat[][SNR 20, TV reg, 64 eval\label{subfig:SNR20I64TVDiff}]{\fbox{\includegraphics[height=\heightNoise\textwidth]{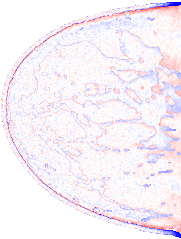}}}\\
\subfloat[][SNR 10, 32 eval\label{subfig:SNR10I32Diff}]{\fbox{\includegraphics[height=\heightNoise\textwidth]{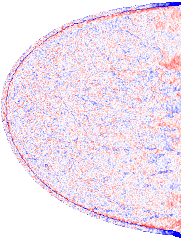}}}
\hfill
\subfloat[][SNR 10, 64 eval\label{subfig:SNR10I64Diff}]{\fbox{\includegraphics[height=\heightNoise\textwidth]{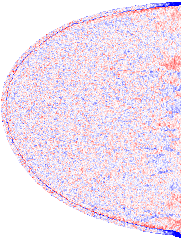}}}
\hfill
\subfloat[][SNR 10, TV reg, 32 eval\label{subfig:SNR10I32TVDiff}]{\fbox{\includegraphics[height=\heightNoise\textwidth]{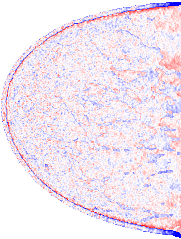}}}
\hfill
\subfloat[][SNR 10, TV reg, 64 eval\label{subfig:SNR10I64TVDiff}]{\fbox{\includegraphics[height=\heightNoise\textwidth]{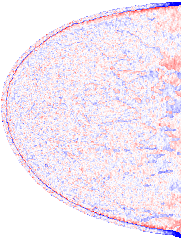}}}\\
\subfloat[][SNR 5, 32 eval\label{subfig:SNR5I32Diff}]{\fbox{\includegraphics[height=\heightNoise\textwidth]{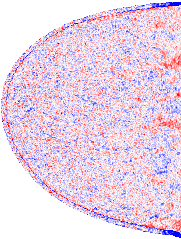}}}
\hfill
\subfloat[][SNR 5, 64 eval\label{subfig:SNR5I64Diff}]{\fbox{\includegraphics[height=\heightNoise\textwidth]{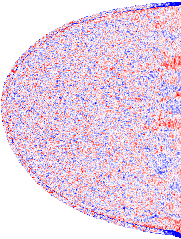}}}
\hfill
\subfloat[][SNR 5, TV reg, 32 eval\label{subfig:SNR5I32TVDiff}]{\fbox{\includegraphics[height=\heightNoise\textwidth]{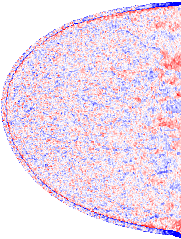}}}
\hfill
\subfloat[][SNR 5, TV reg, 64 eval\label{subfig:SNR5I64TVDiff}]{\fbox{\includegraphics[height=\heightNoise\textwidth]{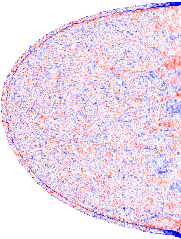}}}\\
\caption{Error plots corresponding to Figure \ref{fig:NoiseStudy}.}
   \label{fig:NoiseStudyDiff}
\end{figure}

\clearpage

\bibliographystyle{plain}
\bibliography{all}

\end{document}